\documentclass[a4paper,11pt,reqno]{amsart}

\usepackage[T1]{fontenc}
\usepackage{lmodern}
\usepackage{a4wide}

\usepackage{amsmath,amssymb,mathtools}
\allowdisplaybreaks[4]
\usepackage{mathrsfs}

\usepackage{enumitem}
\usepackage{needspace}

\usepackage{xcolor}
\usepackage[colorlinks=true,citecolor=teal,linkcolor=blue,urlcolor=blue]{hyperref}
\usepackage{orcidlink}

\makeatletter
\def\@maketitle{\normalfont\normalsize
  \let\@makefnmark\relax  \let\@thefnmark\relax
  \ifx\@empty\@date\else \@footnotetext{\@setdate}\fi
  \ifx\@empty\@subjclass\else \@footnotetext{\@setsubjclass}\fi
  \ifx\@empty\@keywords\else \@footnotetext{\@setkeywords}\fi
  \ifx\@empty\thankses\else \@footnotetext{\def\par{\let\par\@par}\@setthanks}\fi
  \global\topskip42\p@
  \@settitle
  \ifx\@empty\authors \else \@setauthors \fi
  \@setabstract
  \normalsize
  \dimen@34\p@ \advance\dimen@-\baselineskip
  \vskip\dimen@\relax
}
\def\@citestyle{\m@th\upshape\mdseries\color{teal}}
\def\@biblabel#1{\textcolor{teal}{\textbf{[#1]}}}
\makeatother

\numberwithin{equation}{section}

\theoremstyle{plain}

\newtheorem{theorem}{Theorem}[section]
\newtheorem{lemma}[theorem]{Lemma}
\newtheorem{proposition}[theorem]{Proposition}
\newtheorem{corollary}[theorem]{Corollary}

\theoremstyle{remark}

\newtheorem{remark}[theorem]{Remark}

\newcommand{\R}{\mathbb{R}}
\newcommand{\Sph}{\mathbb{S}}

\newcommand{\Ric}{\operatorname{Ric}}
\newcommand{\Rm}{\operatorname{Rm}}
\newcommand{\divg}{\operatorname{div}}
\newcommand{\trg}{\operatorname{tr}}
\newcommand{\Id}{\operatorname{Id}}
\newcommand{\II}{\mathrm{II}}

\newcommand{\inner}[2]{\left\langle #1,#2\right\rangle}
\newcommand{\dmu}{\,\mathrm{d}\mu}

\newcommand{\papertitle}{Rigidity of shrinking Ricci solitons with low-dimensional focal sets}
\newcommand{\papershorttitle}{Shrinking Ricci solitons with low-dimensional focal sets}

\title[\papershorttitle]{\papertitle}

\author[M. A. R. M. Hor\'acio]{Matheus Andrade Ribeiro de Moura Hor\'acio\,\orcidlink{0009-0001-5658-7467}}
\thanks{ORCID iD: \url{https://orcid.org/0009-0001-5658-7467}.}

\address{
Departamento de Matem\'atica,
Universidade de Bras\'ilia,
70910-900, Bras\'ilia, DF,
Brazil
}

\email{matheus.andrade5488@gmail.com}

\date{}



\hypersetup{
  pdftitle={\papertitle},
  pdfauthor={Matheus Andrade Ribeiro de Moura Horacio},
  pdfsubject={Differential geometry},
  pdfkeywords={Gradient Ricci solitons, rigidity, constant scalar curvature, minimum set}
}

\begin{document}
\raggedbottom

\begin{abstract}
We prove that, in every dimension, the focal set of the potential function of a complete gradient shrinking Ricci soliton with constant scalar curvature is totally geodesic whenever it has dimension two or three, and that every complete five-dimensional gradient shrinking Ricci soliton with constant scalar curvature is rigid. The first result follows from a Bochner identity for the covariant derivative of the Ricci tensor at the focal set. For the second, when the scalar curvature is twice the soliton constant, we prove that total geodesy implies rigidity through an integral identity on the level hypersurfaces; the case of scalar curvature three times the soliton constant is treated by a separate argument. As further applications, we combine our methods with the rigidity theorem in the recent preprint of Li, Li, Ou, Qu, Wu and Zhou to obtain rigidity in arbitrary dimension when the scalar curvature is twice or three times the soliton constant, and under integral Weyl pinching for four-dimensional focal sets.
\end{abstract}

\maketitle

\section{Introduction and main results}\label{sec_intro}

A gradient shrinking Ricci soliton $(M^n,g,f,\lambda)$ is a Riemannian manifold $(M^n,g)$ together with a smooth function $f$ and a constant $\lambda>0$ satisfying
\begin{align}\label{fundeq}
	\Ric+\nabla^2f=\lambda g.
\end{align}
Following Petersen and Wylie \cite{pewy}, the soliton is called rigid if it is isometric to a quotient $\bar N\times_\Gamma\R^{n-k}$ of a product of an Einstein manifold $\bar N^k$, with $\Ric_{\bar N}=\lambda g_{\bar N}$, and a Euclidean space, where $f=\frac{\lambda}{2}|x|^2$ on the Euclidean factor, up to a constant, and $\Gamma$ acts freely on $\bar N$ and by orthogonal transformations on $\R^{n-k}$. By \cite[Theorem 1.2]{pewy}, a complete gradient shrinking Ricci soliton is rigid if and only if it has constant scalar curvature and is radially flat, that is, $\Rm(\nabla f,X,\nabla f,X)=0$ for every $X$. Rigidity is also known under several other curvature conditions, see for instance \cite{flgr11,ms13,pw10}.

Cao conjectured that every complete gradient shrinking Ricci soliton with constant scalar curvature is rigid; see \cite{chengzhou}. Fern\'andez-L\'opez and Garc\'ia-R\'io proved that the scalar curvature of such a soliton belongs to $\{0,\lambda,\ldots,n\lambda\}$ \cite[Theorem 1]{fr16}, that $R=\lambda$ does not occur, and that $R=(n-1)\lambda$ gives rigidity \cite[Theorem 10]{fr16}, while $R=0$ and $R=n\lambda$ correspond to the Gaussian soliton and to compact Einstein manifolds \cite[Propositions 3.1 and 3.3]{pewy}. They also proved that such a soliton is rigid as soon as its Ricci tensor has constant rank \cite[Theorem 2]{fr16}. The conjecture was proved in dimension four by Cheng and Zhou \cite{chengzhou}, and another proof of the case $R=2\lambda$ in dimension four was given by Ou, Qu and Wu \cite{ouquwu}; see also \cite{wwres}.
 In dimension five, the remaining values are $R=2\lambda$ and $R=3\lambda$. The case $R=3\lambda$ was considered by Li, Ou, Qu and Wu in the preprint \cite{loqw1}, and the case $R=2\lambda$ was studied by the same authors in \cite{loqw2} under the assumption that the curvature is bounded. The conjecture also holds for K\"ahler solitons in every dimension \cite{hloqwz,horkahler,lowz}, and we refer to \cite{lqw,tranyamazaki} for related results.

Let $(M^n,g,f,\lambda)$ be a complete gradient shrinking Ricci soliton with constant scalar curvature $R=k\lambda$ and nonconstant potential function. After adding a constant to $f$, we have $|\nabla f|^2=2\lambda f$, and the minimum set
\[
	N=f^{-1}(0)
\]
is a compact connected submanifold of dimension $k$ \cite{chengzhou,fr16}, which is the focal set of the transnormal function $f$ in the sense of Wang \cite{wang}. We write $q=n-k$ for its codimension and $\II$ for its second fundamental form. On a rigid soliton $\bar N\times_\Gamma\R^{n-k}$, the minimum set of $f$ is the image of $\bar N\times\{0\}$, naturally identified with $\bar N/\Gamma$, and it is totally geodesic. Our structural result is the following.

\begin{theorem}\label{thm_focal}
Let $(M^n,g,f,\lambda)$ be a complete gradient shrinking Ricci soliton with constant scalar curvature $R=k\lambda$ and nonconstant potential function, and let $N$ be the minimum set of $f$. If $k=2$ or $k=3$, then $N$ is totally geodesic.
\end{theorem}

Its proof does not use any previous rigidity result. Some of its cases also follow from published classifications: the case $q=1$ from the rigidity of solitons with $R=(n-1)\lambda$ \cite[Theorem 10]{fr16}, and the case $(n,R)=(4,2\lambda)$ from the theorem of \cite{chengzhou}. The case $(n,R)=(5,3\lambda)$ was considered in \cite{loqw1}. In the other cases, the conclusion was known only under additional assumptions, such as a K\"ahler structure \cite{hloqwz,horkahler,lowz}. In dimension five, when $R=2\lambda$, we prove that a totally geodesic minimum set forces rigidity by a global argument on the level hypersurfaces of $f$. Together with a different argument for $R=3\lambda$, this gives the following result.

\begin{theorem}\label{thm_five}
Every complete connected five-dimensional gradient shrinking Ricci soliton with constant scalar curvature is rigid. More precisely, with $R=k\lambda$, its universal cover is $\R^5$, $\Sph^2(1/\sqrt{\lambda})\times\R^3$, $\Sph^3(\sqrt{2/\lambda})\times\R^2$, $\widetilde N^4\times\R$ or a compact Einstein manifold according as $k=0$, $2$, $3$, $4$ or $5$, where $\Sph^k(\rho)$ is the round sphere of radius $\rho$ and $\widetilde N^4$ is a compact simply connected Einstein manifold with $\Ric_{\widetilde N}=\lambda g_{\widetilde N}$. The manifold $M$ is the quotient of its universal cover by a finite group, of order at most two when $k=2$.
\end{theorem}

For $R=2\lambda$, Theorem \ref{thm_five} follows from Theorem \ref{thm_focal} and the argument on the level hypersurfaces in Section \ref{sec_level}; for $R=3\lambda$, Section \ref{sec_three} gives an argument independent of Theorem \ref{thm_focal}. Neither argument uses \cite{lloqwz}, \cite{loqw1} or \cite{loqw2}. The bound on the order of the quotient is attained by $\mathbb{RP}^2\times\R^3$, and the example $\mathbb{CP}^2\times\R$ shows that the Einstein factor cannot be specified further when $k=4$.

\subsection*{Further applications}

While this article was being prepared, Li, Li, Ou, Qu, Wu and Zhou posted the preprint \cite{lloqwz}, which proves that, in every dimension, a complete gradient shrinking Ricci soliton with constant scalar curvature is rigid whenever the minimum set of its potential function is totally geodesic \cite[Corollary 1.1]{lloqwz}. Theorems \ref{thm_focal} and \ref{thm_five} do not use \cite{lloqwz}, which we use only to pass from total geodesy to rigidity in Corollary \ref{cor_alldim} and Theorems \ref{thm_four} and \ref{thm_harmonic_weyl}. Combined with Theorem \ref{thm_focal}, the theorem of \cite{lloqwz} gives the following consequence.

\begin{corollary}\label{cor_alldim}
Let $(M^n,g,f,\lambda)$ be a complete connected gradient shrinking Ricci soliton with constant scalar curvature $R\in\{2\lambda,3\lambda\}$. Then it is rigid. More precisely, if $R=k\lambda$, its universal cover is $\Sph^k(\sqrt{(k-1)/\lambda})\times\R^{n-k}$, and $M$ is the quotient of this product by a finite group of isometries which preserve the potential function. Consequently, every complete gradient shrinking Ricci soliton with constant scalar curvature $R\leq3\lambda$ is rigid, and in dimension six the conjecture of Cao holds, except possibly when $R=4\lambda$.
\end{corollary}

In the K\"ahler case, the minimum set is the zero set of the Killing field $J\nabla f$, hence totally geodesic, so the theorem of \cite{lloqwz} also connects with the rigidity of K\"ahler shrinkers with constant scalar curvature \cite{hloqwz,horkahler,lowz}; see Remark \ref{rem_kahler}.

When $R=4\lambda$, the minimum set $N$ is a closed minimal four-manifold, and the integral identity behind Theorem \ref{thm_focal} contains a term involving the Weyl tensor $W^N$ of the metric induced on $N$, which has no sign. For an orthonormal frame $\{e_i\}$ of $TN$, we write
\[
	W^N(\II,\II)=\sum_{i,j,\ell,m}W^N(e_i,e_j,e_\ell,e_m)\inner{\II(e_i,e_\ell)}{\II(e_j,e_m)},
\]
and we denote by $D$ the connection induced on $T^*N\otimes T^*N\otimes E$, where $E$ is the normal bundle of $N$, by the Levi-Civita connection of $N$ and the normal connection. The following result controls this term by an integral pinching condition.

\begin{theorem}\label{thm_four}
Let $(M^n,g,f,\lambda)$ be a complete gradient shrinking Ricci soliton with constant scalar curvature $R=4\lambda$ and nonconstant potential function, and let $N$ be the minimum set of $f$. Then
\begin{align}\label{eq_four_estimate}
	11\int_NW^N(\II,\II)\dmu_N\geq\int_N\Big(|D\II|^2+\frac23\lambda|\II|^2+\frac{1}{72}|\II|^4\Big)\dmu_N,
\end{align}
and, if $n=6$, the coefficient $\frac{1}{72}$ in \eqref{eq_four_estimate} can be replaced by $\frac{11}{8}$. Consequently, if
\begin{align}\label{eq_four_pinching}
	\int_NW^N(\II,\II)\dmu_N\leq\frac{2\lambda}{33}\int_N|\II|^2\dmu_N,
\end{align}
then $N$ is totally geodesic and $M$ is rigid.
\end{theorem}

In particular, $M$ is rigid whenever $\int_NW^N(\II,\II)\dmu_N\leq0$. For the second result, we regard $W^N$ as a self-adjoint endomorphism of $\Lambda^2TN$ by
\[
	\inner{W^N(e_i\wedge e_j)}{e_\ell\wedge e_m}=W^N(e_i,e_j,e_\ell,e_m),
\]
where $\{e_i\wedge e_j\}_{i<j}$ is orthonormal, as for the level sets in Section \ref{sec_level}. We write $\|W^N\|$ and $\|\nabla^NW^N\|$ for the induced Hilbert--Schmidt norms, so that $\|W^N\|^2=\frac14|W^N|^2$, $\trg_{\Lambda^2}((W^N)^3)$ for the trace of the cube of this endomorphism, and $\delta_NW^N$ for the divergence of $W^N$ on $(N,g_N)$.

\begin{theorem}\label{thm_harmonic_weyl}
Let $(M^n,g,f,\lambda)$ be a complete gradient shrinking Ricci soliton with constant scalar curvature $R=4\lambda$ and nonconstant potential function, and let $N$ be the minimum set of $f$. If $\delta_NW^N=0$ and
\begin{align}\label{eq_weyl_condition}
	6\int_N\trg_{\Lambda^2}\big((W^N)^3\big)\dmu_N\leq\int_N\|\nabla^NW^N\|^2\dmu_N+2\lambda\int_N\|W^N\|^2\dmu_N,
\end{align}
then $N$ is totally geodesic and $M$ is rigid.
\end{theorem}

The additional curvature conditions of Theorems \ref{thm_four} and \ref{thm_harmonic_weyl} hold on every rigid soliton with $R=4\lambda$. When $n=5$, the value $R=4\lambda$ already gives rigidity \cite[Theorem 10]{fr16}, so the first new case is $n=6$, which is not covered by Theorem \ref{thm_focal} and Corollary \ref{cor_alldim}.

\subsection*{Strategy of proof}

The proof of Theorem \ref{thm_focal} takes place at the minimum set. At the points of $N$, we have $\Ric=\lambda P_{TN}$ and $|\nabla\Ric|^2=3\lambda^2|\II|^2$ (Lemma \ref{lem_first_order}), so a Bochner formula for $\nabla\Ric$ at $N$ is a natural way to study total geodesy. We obtain it by applying the Laplacian to the Bochner formula \eqref{eq_bochner_global} for $\Ric$, which holds on $M$ since $R$ and $|\Ric|^2$ are constant. Separating tangent and normal entries, using the Gauss, Codazzi and Ricci equations, and integrating over the closed manifold $N$, we obtain $\int_NJ\dmu_N=0$, in every dimension and codimension (Theorem \ref{thm_general}). After completing the squares in the curvature components of $M$ with two normal entries, $J$ consists of squares, a positive multiple of $\lambda|\II|^2$, and terms involving the intrinsic curvature of $N$ and quartic expressions in the shape operators.

When $k=2$, the shape operators are trace free symmetric endomorphisms of a plane, so that $S_\alpha S_\beta+S_\beta S_\alpha$ is a multiple of the identity, and the remaining terms are nonnegative (Theorem \ref{thm_two}). When $k\geq3$, a Weitzenb\"ock formula for $\II$ gives a second vanishing integral $\int_NI\dmu_N=0$ (Lemma \ref{lem_simons}). When $k=3$, subtracting four times this identity eliminates $|D\II|^2$ and leaves $\lambda|\II|^2$ with coefficient one, and a polarized Cayley--Hamilton identity makes the remaining quartic terms nonnegative (Theorem \ref{thm_three}). In both cases, the term $\lambda|\II|^2$ forces $\II=0$. When $k\geq4$, $J$ contains the pairing $W^N(\II,\II)$ of the Weyl tensor of $N$ with $\II$, which has no sign. Constant scalar curvature does not force $W^N=0$, even on rigid solitons, but the pairing vanishes on every rigid soliton (Remark \ref{rem_four}), and Theorems \ref{thm_four} and \ref{thm_harmonic_weyl} control it.

Theorem \ref{thm_five} uses Theorem \ref{thm_focal} only when $R=2\lambda$, where $N$ is a surface of codimension three, and Section \ref{sec_level} proves that total geodesy implies rigidity in this case. The level hypersurfaces $\Sigma_r$ of $f$ are closed four-manifolds, and the Chern--Gauss--Bonnet formula on $\Sigma_r$, together with the Gauss--Bonnet theorem on $N$, gives the integral of the squared norm of the Weyl tensor $W_\Sigma$ of $\Sigma_r$ in terms of $\int_N|\II|^2\dmu_N$ and of the normal derivative of $\Ric$. The Bochner formula restricted to $\Sigma_r$ and an integral identity for divergence free tensors give the integral of the pairing of $W_\Sigma$ with the trace free part $B$ of the restriction of $\Ric$ to $T\Sigma_r$, which is represented by the Weyl component of $B\wedge B$. Expanding the square of a suitable combination of $W_\Sigma$ and this component gives the identity of Theorem \ref{thm_identity}, which expresses $8\pi\lambda^2\int_N|\II|^2\dmu_N$ as the integral over $\Sigma_r$ of four nonnegative terms and sharpens the Chern--Gauss--Bonnet estimates used in \cite{loqw1,loqw2}. When $N$ is totally geodesic, these terms vanish, which forces the shape operator $S$ of the level sets to satisfy $S'+S^2=0$. Hence radial flatness and \cite[Theorem 1.2]{pewy} give rigidity (Theorem \ref{thm_level}).

For $R=3\lambda$, Section \ref{sec_three} uses an argument independent of Theorem \ref{thm_focal}. The rescaled shape operator $U=rS$ of the level sets satisfies $\trg U=|U|^2=1$, so $1-\trg U^3$ is a nonnegative cubic defect which vanishes exactly when $U$ is a rank one orthogonal projection. The Chern--Gauss--Bonnet formula on $\Sigma_r$ and an algebraic identity for this defect show that the derivative of a weighted average of the defect over $\Sigma_r$ is a sum of squares \eqref{eq_W_average}. Since this weighted average tends to zero at infinity, the defect vanishes. Then $U$ is a rank one orthogonal projection, and the same Riccati computation gives radial flatness, so that radial flatness and \cite[Theorem 1.2]{pewy} again give rigidity.

Sections \ref{sec_prelim} to \ref{sec_bochner} derive the Bochner identity at $N$ and its integral form, Section \ref{sec_lowdim} proves Theorem \ref{thm_focal} and Corollary \ref{cor_alldim}, Sections \ref{sec_level} and \ref{sec_three} prove Theorem \ref{thm_five}, and Section \ref{sec_four} proves Theorems \ref{thm_four} and \ref{thm_harmonic_weyl}. Appendix \ref{app_conn} contains the computation of the connection terms in Lemma \ref{lem_X}, and Appendix \ref{app_notation} lists the notation of Sections \ref{sec_prelim} to \ref{sec_lowdim} and \ref{sec_four}.

\subsection*{Acknowledgements}

This study was financed by the Coordena\c{c}\~ao de Aperfei\c{c}oamento de Pessoal de N\'ivel Superior -- Brasil (CAPES) -- Finance Code 001. The author thanks the Departamento de Matem\'atica of the Universidade de Bras\'ilia, where this work was carried out. He is particularly grateful to his doctoral advisor, Jo\~ao Paulo dos Santos, for his valuable guidance and support. He also thanks Detang Zhou for helpful discussions on the conjecture of Cao.

\subsection*{AI Disclosure}

Publicly available large language models were used in preparing this manuscript, and their contributions were limited exclusively to the following: symbolic and algebraic manipulations, namely finding the coefficients that express the integral identities as sums of squares, computer-algebra verification of the index computations, and the application of standard linear algebra identities in Section \ref{sec_lowdim}, such as the Cayley--Hamilton identity in dimension three; adapting, at the author's request, the monotonicity argument of \cite[Corollary 1.6]{horlevel} to the proof of Theorem \ref{thm_three_five}; and extensive literature searches and proofreading. The author wrote the entire text and every mathematical argument and proof in it, checked all of them, and takes full responsibility for all mathematical arguments and the final text.

\section{Preliminaries}\label{sec_prelim}

All manifolds are assumed to be connected. We adopt the convention
\[
	\Rm(X,Y,Z,T)=\inner{\nabla_{Y}\nabla_{X}Z-\nabla_{X}\nabla_{Y}Z+\nabla_{[X,Y]}Z}{T},
\]
so that $\Rm(X,Y,X,Y)$ is the sectional curvature of $\mathrm{span}\{X,Y\}$ for orthonormal $X$ and $Y$, and $\Ric(Y,Z)=\sum_a\Rm(e_a,Y,e_a,Z)$ for an orthonormal frame $\{e_a\}$. With $R_{abcd}=\Rm(e_a,e_b,e_c,e_d)$, the Ricci identity reads
\begin{align}\label{eq_ricci_identity}
	\nabla^2_{a,b}\phi_c-\nabla^2_{b,a}\phi_c=\sum_dR_{abcd}\phi_d
\end{align}
for a one-form $\phi$, and similarly for each entry of a tensor, and the contracted second Bianchi identity reads $\sum_a\nabla_aR_{abcd}=\nabla_c\Ric_{bd}-\nabla_d\Ric_{bc}$. We use the same symbol for a symmetric bilinear form and the corresponding self-adjoint endomorphism, and, unless otherwise stated, the norms of tensors are the full tensor norms. For symmetric $A$ and $B$ we write
\[
	\Rm(A)_{bd}=\sum_{a,c}R_{abcd}A_{ac},\qquad \Rm(A,B)=\sum_{a,b,c,d}R_{abcd}A_{ac}B_{bd},
\]
so that $\Rm(A,A)=\sum_{a,b}R_{abab}\alpha_a\alpha_b$ if $\{e_a\}$ is an orthonormal basis of eigenvectors of $A$ with eigenvalues $\alpha_a$. Recall that a gradient Ricci soliton satisfies
\begin{align}\label{eq_identities}
	R+\Delta f=n\lambda,\qquad \nabla R=2\Ric(\nabla f,\cdot),\qquad R+|\nabla f|^2-2\lambda f=\mathrm{const},
\end{align}
see \cite{pewy}, as well as
\begin{align}\label{eq_soliton_commutation}
\begin{split}
	(\nabla_X\Ric)(Y,Z)-(\nabla_Y\Ric)(X,Z)&=\Rm(X,Y,\nabla f,Z),\\
	\Delta_f\Ric&=2\lambda\Ric-2\Rm(\Ric),
\end{split}
\end{align}
where $\Delta_f=\Delta-\nabla_{\nabla f}$, see \cite[Proposition 2.1]{elm} and \cite{pewy}.

Throughout Sections \ref{sec_prelim} to \ref{sec_lowdim} and in Section \ref{sec_four}, $(M^n,g,f,\lambda)$ is a complete gradient shrinking Ricci soliton with constant scalar curvature $R=k\lambda$ and nonconstant potential function. By \cite[Theorems 1 and 10]{fr16}, $0\leq k\leq n-1$ and $k\neq1$; if $k=0$, then $\Ric=0$ and $M$ is the Gaussian soliton \cite[Propositions 3.1 and 3.3]{pewy}, so we assume $2\leq k\leq n-1$. By \eqref{eq_identities}, $\Ric(\nabla f,\cdot)=0$, and after adding a constant to $f$ we have $|\nabla f|^2=2\lambda f$. Hence
\begin{align}\label{eq_laplacian_f}
	\Delta f=q\lambda,\qquad f\geq0,\qquad q=n-k.
\end{align}
Tracing the second identity in \eqref{eq_soliton_commutation} gives $\Delta_fR=2\lambda R-2|\Ric|^2$. Since $R$ is constant, $|\Ric|^2=\lambda R$ is constant as well, and the same identity gives
\begin{align}\label{eq_bochner_global}
	0=\frac12\Delta_f|\Ric|^2=|\nabla\Ric|^2+2\lambda|\Ric|^2-2\Rm(\Ric,\Ric).
\end{align}

\begin{proposition}\label{prop_structure}
The set $N=f^{-1}(0)$ is a compact connected embedded submanifold of $M$ of dimension $k$, $f$ attains its minimum exactly on $N$, and the function $r=\sqrt{2f/\lambda}$ is the distance to $N$. At the points of $N$, $\nabla^2f=\lambda P_E$, so that $\Ric=\lambda P_{TN}$ by \eqref{fundeq}, where $P_E$ and $P_{TN}$ are the orthogonal projections onto the normal bundle $E$ of $N$ and onto $TN$, and $N$ is minimal. Moreover, the normal exponential map $\exp\colon E\to M$ is a diffeomorphism, and, if $\mathrm{d}\mu_N$ denotes the volume measure of $N$ and $z$ the components of a normal vector in a local orthonormal frame of $E$, then the volume measure of $M$, pulled back by $\exp$, is $\mathrm{d}\mu_N\,\mathrm{d}z$.
\end{proposition}

\begin{proof}
Since $f$ is nonconstant, $M$ is noncompact: otherwise $\nabla f$ would vanish at a maximum point of $f$, and $|\nabla f|^2=2\lambda f$ would force $f\equiv0$. After scaling $\lambda$ to $\frac12$, \cite[Theorem 2.2]{chengzhou} shows that $N$ is a compact connected minimal submanifold of dimension $k$, that $r$ is the distance to $N$, and that $\nabla^2f=\lambda P_E$ at the points of $N$; see also \cite{fr16} and \cite{wang}. Since $N$ is compact, its inclusion is an embedding, and, since $f\geq0$, $f$ attains its minimum exactly on $N$.
On $M\setminus N$ we have $|\nabla r|=1$, so the integral curves of $\nabla r$ are unit speed geodesics, which near $N$ are the normal geodesics of $N$. By completeness, the integral curve through $x\in M\setminus N$ reaches $N$ backward after length $r(x)$ and is defined forward for every positive value of $r$. Hence $\exp$ is a bijection from $E$ onto $M$, and it is a diffeomorphism, since its restriction to $\{|z|=r\}$ is the composition of the normal exponential map on a small tube with the flow of $\nabla r$.

Let $\theta$ be the density of the pulled back volume measure with respect to $\mathrm{d}\mu_N\,\mathrm{d}z$, which does not depend on the choice of the orthonormal frame of $E$. By the Gauss lemma, $r\nabla r=z^\alpha\partial_{z^\alpha}$, so that $\nabla f=\lambda z^\alpha\partial_{z^\alpha}$ and
\[
	q\lambda=\Delta f=\lambda\left(q+z^\alpha\partial_{z^\alpha}\log\theta\right).
\]
Thus $\theta$ is constant along the rays $t\mapsto tz$, and $\theta=1$, since $\theta=1$ on $N$.
\end{proof}

We fix $p\in N$, an orthonormal basis $\{e_i\}$ of $T_pN$ and an orthonormal basis $\{e_\alpha\}$ of $E_p$. Indices $i,j,\ell,m$ are tangent, indices $\alpha,\beta,\gamma,\delta$ are normal, and indices $a,b,c,d,e,s,x$ run over all directions. Repeated indices are summed. We write $\II(X,Y)=(\nabla_XY)^\perp$ for the second fundamental form of $N$, $h^\alpha_{ij}=\inner{\II(e_i,e_j)}{e_\alpha}$, $S_\alpha$ for the shape operator in the direction of $e_\alpha$, so that $\inner{S_\alpha X}{Y}=h^\alpha(X,Y)$, $D$ for the connection induced on $T^*N\otimes T^*N\otimes E$ by the Levi-Civita connection of $N$ and the normal connection $\nabla^\perp$, and $R^\perp(X,Y)=\nabla^\perp_X\nabla^\perp_Y-\nabla^\perp_Y\nabla^\perp_X-\nabla^\perp_{[X,Y]}$ for the curvature of the normal connection. Products of shape operators are written as matrices in the basis $\{e_i\}$, so that $(S_\alpha S_\beta)_{ij}=\sum_\ell h^\alpha_{i\ell}h^\beta_{\ell j}$. We also use the following quantities at $p$:
\begin{center}
\renewcommand{\arraystretch}{1.6}
\begin{tabular}{@{}l@{\qquad}l@{}}
	shape operators & $t=|\II|^2$,\quad $G_{\alpha\beta}=\trg(S_\alpha S_\beta)$,\quad $c=|G|^2$,\quad $T=\sum_\alpha S_\alpha^2$,\\
	& $P_{\alpha\beta}=S_\alpha S_\beta+S_\beta S_\alpha$,\quad $K_{\alpha\beta}=S_\alpha S_\beta-S_\beta S_\alpha$,\\
	& $\sigma_{ij\ell m}=\inner{\II_{ij}}{\II_{\ell m}}+\inner{\II_{i\ell}}{\II_{jm}}+\inner{\II_{im}}{\II_{j\ell}}$,\\
	divergence of $\II$ & $d^\alpha_j=\sum_i(D_{e_i}h^\alpha)_{ij}$,\\
	curvature of $M$ & $L_{\alpha\beta ij}=R_{\alpha i\beta j}$,\quad $V_{\alpha\beta ij}=R_{\alpha\beta ij}$,\quad $Y_{\alpha i\beta\gamma}=R_{\alpha i\beta\gamma}$,\\
	& $\rho_{ij}=\sum_\ell R_{\ell i\ell j}$,\quad $\tilde\rho=\Rm(P_{TN})$, so that $\tilde\rho_{bd}=\sum_iR_{ibid}$.
\end{tabular}
\end{center}
The tensor $\sigma$ is totally symmetric on $T_pN$. Let $L^s$ be the symmetrization of $L$ in $(i,j)$ and $L^t_{\alpha\beta ij}=L_{\alpha\beta ji}$. By the first Bianchi identity, $L_{\alpha\beta ij}-L_{\alpha\beta ji}=V_{\alpha\beta ij}$, so that
\begin{align}\label{eq_N_norms}
	|L|^2=|L^s|^2+\tfrac14|V|^2,\qquad \inner{L}{L^t}=|L^s|^2-\tfrac14|V|^2.
\end{align}
Finally, pairings such as $\inner{K}{V}=\sum(K_{\alpha\beta})_{ij}V_{\alpha\beta ij}$ and $\inner{P}{L^s}=\sum(P_{\alpha\beta})_{ij}L^s_{\alpha\beta ij}$ are full contractions over all indices. Further quantities are introduced where they are first used, and Appendix \ref{app_notation} lists all of them, with the place of their definition.

\needspace{4\baselineskip}
\begin{lemma}\label{lem_first_order}
At $p\in N$:
\begin{enumerate}[label=\textup{(\roman*)},leftmargin=*]
\item $\Ric=\lambda P_{TN}$, and $\nabla\Ric$ is totally symmetric; its only nonzero components are those with exactly one normal entry, and $(\nabla\Ric)(e_i,e_j,e_\alpha)=\lambda h^\alpha_{ij}$. In particular, $|\nabla\Ric|^2=3\lambda^2t$.
\item (Codazzi equation) $R_{i\ell j\alpha}=(D_\ell h^\alpha)_{ij}-(D_ih^\alpha)_{\ell j}$, and $\sum_iR_{i\alpha ij}=-d^\alpha_j$.
\item (Ricci equation) $\inner{R^\perp(e_i,e_j)e_\beta}{e_\alpha}=V_{\alpha\beta ij}+(K_{\alpha\beta})_{ij}$.
\item For all $a,b,c,e$, $\nabla^2_{e,a}\Ric_{bc}-\nabla^2_{e,b}\Ric_{ac}=\lambda\Rm(e_a,e_b,P_Ee_e,e_c)$.
\end{enumerate}
\end{lemma}

\begin{proof}
(i) The first assertion is Proposition \ref{prop_structure}. Since $\nabla f=0$ at $p$, the first identity in \eqref{eq_soliton_commutation} shows that $\nabla\Ric$ is symmetric in its first two entries, hence totally symmetric. For tangent $X,Y,Z$ and normal $\nu,\eta$, differentiating $\Ric=\lambda P_{TN}$ along $N$ gives
\[
	(\nabla_X\Ric)(Y,Z)=0,\qquad (\nabla_X\Ric)(\nu,\eta)=0,\qquad (\nabla_X\Ric)(Y,\nu)=\lambda\inner{\II(X,Y)}{\nu},
\]
since $(\nabla_X\nu)^{TN}=-S_\nu X$. By total symmetry, this determines every component with at least one tangent entry. For a unit normal $\nu$ at $p$, let $\gamma(s)=\exp_p(s\nu)$. By Proposition \ref{prop_structure}, $\nabla f=\lambda s\gamma'$ along $\gamma$, so $\Ric(\gamma',\gamma')=0$ for $s>0$. Differentiating at $s=0$ gives $(\nabla_\nu\Ric)(\nu,\nu)=0$, and polarization of this totally symmetric cubic form on $E_p$ shows that the components with three normal entries vanish. Counting the three positions of the normal entry gives $|\nabla\Ric|^2=3\lambda^2t$.

(ii) The Codazzi equation is the normal part of the definition of $\Rm$. Summing over $i=j$ and using the minimality of $N$ (Proposition \ref{prop_structure}) gives $\sum_iR_{i\ell i\alpha}=-d^\alpha_\ell$, and $R_{i\alpha i\ell}=R_{i\ell i\alpha}$.

(iii) Taking the normal part of the definition of $\Rm$, with $\nabla_X\xi=\nabla^\perp_X\xi-S_\xi X$ for a normal field $\xi$, we obtain
\[
	\Rm(X,Y,\xi,\eta)=-\inner{R^\perp(X,Y)\xi}{\eta}+\inner{[S_\xi,S_\eta]X}{Y}.
\]
We take $X=e_i$, $Y=e_j$, $\xi=e_\alpha$ and $\eta=e_\beta$. In the matrix convention fixed above, $\inner{[S_\alpha,S_\beta]e_i}{e_j}=-(K_{\alpha\beta})_{ij}$, and (iii) follows.

(iv) Differentiate the first identity in \eqref{eq_soliton_commutation} in the direction $e_e$ and use $\nabla f=0$, $\nabla_{e_e}\nabla f=\lambda P_Ee_e$ at $p$.
\end{proof}

\section{The Hessian of the Ricci tensor at the minimum set}\label{sec_hessian}

The Bochner identity of Section \ref{sec_bochner} contains $|\nabla^2\Ric|^2$ at the points of $N$. In this section, we express every component of $\nabla^2\Ric$ at $p\in N$ with at least one tangent entry through $\II$, $D\II$ and the curvature of $M$ along $N$ (Lemma \ref{lem_hessian}), which gives the norm \eqref{eq_hessian_norm}. Lemma \ref{lem_hess_two} relates the curvature of $M$ along $N$ to the intrinsic curvature of $N$.

Write $\Lambda_{eabc}=(\nabla^2_{e,a}\Ric)(e_b,e_c)$ at $p$. It is symmetric in $(b,c)$, and, by \eqref{eq_ricci_identity} and Lemma \ref{lem_first_order}(i), $\Lambda_{eabc}-\Lambda_{aebc}=\sum_dR_{eabd}\Ric_{dc}+\sum_dR_{eacd}\Ric_{bd}$, that is,
\begin{align}\label{eq_lambda_comm}
	\Lambda_{eabc}-\Lambda_{aebc}=\lambda\big(R_{eabc}[c\in TN]+R_{eacb}[b\in TN]\big),
\end{align}
where $[c\in TN]$ is one if $e_c$ is tangent and zero otherwise. By Lemma \ref{lem_first_order}(iv), $\Lambda_{eabc}$ is totally symmetric in $(a,b,c)$ when $e$ is tangent. We write $h^\alpha_{ij,\ell}=(D_\ell h^\alpha)_{ij}$ for the components of $D\II$.

\begin{lemma}\label{lem_hessian}
At $p\in N$:
\begin{enumerate}[label=\textup{(\roman*)},leftmargin=*]
\item $\Lambda_{\ell ijm}=-\lambda\sigma_{\ell ijm}$;
\item every component with exactly one normal entry equals $\lambda h^\alpha_{jm,\ell}$, where $\ell$ is the first tangent entry and $(j,m)$ are the other two tangent entries;
\item the components with two normal entries are
\begin{gather*}
	\Lambda_{ij\alpha\beta}=\Lambda_{i\alpha j\beta}=\Lambda_{i\alpha\beta j}=\lambda(P_{\alpha\beta})_{ij},\qquad \Lambda_{\alpha i\beta j}=\Lambda_{\alpha ij\beta}=\lambda(P_{\alpha\beta}+L_{\alpha\beta})_{ij},\\
	\Lambda_{\alpha\beta ij}=\lambda(P_{\alpha\beta}+2L^s_{\alpha\beta})_{ij};
\end{gather*}
\item the components with three normal entries are
\[
	\Lambda_{i\alpha\beta\gamma}=\Lambda_{\alpha i\beta\gamma}=0,\qquad \Lambda_{\alpha\beta i\gamma}=\Lambda_{\alpha\beta\gamma i}=\lambda Y_{\beta i\alpha\gamma}.
\]
\end{enumerate}
Consequently, with $F_{\alpha\beta\gamma\delta}=\Lambda_{\alpha\beta\gamma\delta}$,
\begin{align}\label{eq_hessian_norm}
	\frac{1}{\lambda^2}|\nabla^2\Ric|^2=|\sigma|^2+4|D\II|^2+6|P|^2+8\inner{P}{L^s}+6|L^s|^2+\frac12|V|^2+2|Y|^2+\frac{1}{\lambda^2}|F|^2.
\end{align}
\end{lemma}

\begin{proof}
Extend the frames along $N$ so that $\nabla^Ne_i=0$ and $\nabla^\perp e_\alpha=0$ at $p$. Then $\nabla_{e_\ell}e_i=\II(e_\ell,e_i)$ and $\nabla_{e_\ell}e_\alpha=-S_\alpha e_\ell$ at $p$. By Lemma \ref{lem_first_order}(i), the components of $\nabla\Ric$ with zero, two or three normal entries vanish identically along $N$, and $(\nabla\Ric)(e_\alpha,e_j,e_m)=\lambda h^\alpha_{jm}$. Differentiating these components along $N$ at $p$, we obtain
\begin{align*}
	\Lambda_{\ell ijm}&=-(\nabla_{\II(e_\ell,e_i)}\Ric)(e_j,e_m)-(\nabla_{e_i}\Ric)(\II(e_\ell,e_j),e_m)\\
	&\quad-(\nabla_{e_i}\Ric)(e_j,\II(e_\ell,e_m))=-\lambda\sigma_{\ell ijm},\\
	\Lambda_{\ell\alpha jm}&=e_\ell\big(\lambda h^\alpha_{jm}\big)+(\nabla_{S_\alpha e_\ell}\Ric)(e_j,e_m)-(\nabla_{e_\alpha}\Ric)(\II(e_\ell,e_j),e_m)\\
	&\quad-(\nabla_{e_\alpha}\Ric)(e_j,\II(e_\ell,e_m))=\lambda h^\alpha_{jm,\ell},\\
	\Lambda_{ij\alpha\beta}&=-(\nabla_{\II(e_i,e_j)}\Ric)(e_\alpha,e_\beta)+(\nabla_{e_j}\Ric)(S_\alpha e_i,e_\beta)\\
	&\quad+(\nabla_{e_j}\Ric)(e_\alpha,S_\beta e_i)=\lambda(S_\alpha S_\beta+S_\beta S_\alpha)_{ij},\\
	\Lambda_{i\alpha\beta\gamma}&=(\nabla_{S_\alpha e_i}\Ric)(e_\beta,e_\gamma)+(\nabla_{e_\alpha}\Ric)(S_\beta e_i,e_\gamma)\\
	&\quad+(\nabla_{e_\alpha}\Ric)(e_\beta,S_\gamma e_i)=0,
\end{align*}
where in each line the terms which involve components of $\nabla\Ric$ with zero or two normal entries have been dropped. The total symmetry of $\Lambda_{\ell\cdot\cdot\cdot}$ moves the normal entries among the last three slots, which gives the remaining components with a tangent first entry. The components with a normal first entry follow from \eqref{eq_lambda_comm}: for instance, $\Lambda_{\alpha i\beta j}=\Lambda_{i\alpha\beta j}+\lambda R_{\alpha i\beta j}$, and
\[
	\Lambda_{\alpha \ell jm}=\Lambda_{\ell\alpha jm}+\lambda(R_{\alpha \ell jm}+R_{\alpha \ell mj})=\lambda h^\alpha_{jm,\ell}.
\]
Finally, Lemma \ref{lem_first_order}(iv) with $e=e_\alpha$ gives $\Lambda_{\alpha\beta ij}=\Lambda_{\alpha i\beta j}+\lambda R_{\beta i\alpha j}$ and $\Lambda_{\alpha\beta i\gamma}=\Lambda_{\alpha i\beta\gamma}+\lambda R_{\beta i\alpha\gamma}$, which proves (iii) and (iv).

For the norm, the four components with one normal entry contribute $4|D\II|^2$, the six components with two normal entries contribute $3|P|^2+2|P+L|^2+|P+2L^s|^2$, which equals $6|P|^2+8\inner{P}{L^s}+6|L^s|^2+\frac12|V|^2$ by \eqref{eq_N_norms}, since $P$ is symmetric in $(i,j)$, and the components with three normal entries contribute $2|Y|^2$.
\end{proof}

\begin{lemma}\label{lem_hess_two}
At $p\in N$:
\begin{enumerate}[label=\textup{(\roman*)},leftmargin=*]
\item $\sum_iR_{i\alpha i\beta}=-\frac32G_{\alpha\beta}$ and $\sum_\gamma R_{\gamma\alpha\gamma\beta}=\frac32G_{\alpha\beta}$;
\item $\rho=\Ric_N+T$, $\trg\rho=k\lambda+\frac32t$, and $|\tilde\rho|^2=|\rho|^2+2|d|^2+\frac94c$.
\end{enumerate}
In particular, the scalar curvature of $N$ is $R_N=k\lambda+\frac12t$. If $k=2$, then $\rho=(\lambda+\frac34t)\Id$ and the Gaussian curvature of $N$ is $K_N=\lambda+\frac14t$.
\end{lemma}

\begin{proof}
(i) Since $|\Ric|^2=\lambda R$ is constant, its second derivative in the directions $e_\alpha,e_\beta$ vanishes at $p$:
\[
	0=\sum_{b,c}\Lambda_{\alpha\beta bc}\Ric_{bc}+\sum_{b,c}(\nabla_\alpha\Ric)_{bc}(\nabla_\beta\Ric)_{bc}=\lambda\sum_i\Lambda_{\alpha\beta ii}+\lambda^2G_{\alpha\beta}.
\]
By Lemma \ref{lem_hessian}(iii), $\sum_i\Lambda_{\alpha\beta ii}=\lambda\big(2G_{\alpha\beta}+2\sum_iR_{\alpha i\beta i}\big)$, which gives the first identity. The second follows from $\Ric(e_\alpha,e_\beta)=0$.

(ii) The Gauss equation $R^N_{ij\ell m}=R_{ij\ell m}+\inner{\II_{i\ell}}{\II_{jm}}-\inner{\II_{im}}{\II_{j\ell}}$ and the minimality of $N$ (Proposition \ref{prop_structure}) give $\Ric_N=\rho-T$. Tracing $\Ric(e_i,e_i)=\lambda$ over $i$ and using (i), $k\lambda=\trg\rho-\frac32t$. By Lemma \ref{lem_first_order}(ii) and (i), the blocks of $\tilde\rho$ are
\[
	\tilde\rho_{ij}=\rho_{ij},\qquad \tilde\rho_{\alpha j}=-d^\alpha_j,\qquad \tilde\rho_{\alpha\beta}=-\frac32G_{\alpha\beta},
\]
which gives the formula for $|\tilde\rho|^2$. The scalar curvature of $N$ is $\trg\Ric_N=\trg\rho-t$.

If $k=2$, then $\rho$ is a multiple of the identity on the two-dimensional space $T_pN$, so $\rho=\frac12(\trg\rho)\Id=(\lambda+\frac34t)\Id$, and $K_N=\frac12R_N=\lambda+\frac14t$.
\end{proof}

When $k=2$, the Gauss--Bonnet theorem and Lemma \ref{lem_hess_two} give
\begin{align}\label{eq_gauss_bonnet_N}
	2\pi\chi(N)=\lambda|N|+\frac14\int_N|\II|^2\dmu_N.
\end{align}

\section{A Bochner identity at the minimum set}\label{sec_bochner}

In this section, we prove the integral identity $\int_NJ\dmu_N=0$ of Theorem \ref{thm_general}, on which Section \ref{sec_lowdim} is based, and the companion identity $\int_NI\dmu_N=0$ of Lemma \ref{lem_simons}. We first obtain a pointwise Bochner identity for $\nabla\Ric$ at the points of $N$ (Proposition \ref{prop_bochner}), and then integrate it over $N$.

The identity \eqref{eq_bochner_global} expresses $|\nabla\Ric|^2$ through the curvature at every point of $M$. We apply the Laplacian to it and evaluate the result at $p\in N$. Since $|\Ric|^2$ is constant,
\begin{align}\label{eq_bochner_laplacian}
	|\nabla^2\Ric|^2+\inner{\nabla\Ric}{\Delta\nabla\Ric}=(\Delta\Rm)(\Ric,\Ric)+2\Rm(\Delta\Ric,\Ric)+2\Theta+4\lambda^2\Omega,
\end{align}
where $\Delta\Rm$ is the rough Laplacian of the curvature tensor,
\[
	\Theta=\sum_e\Rm(\nabla_{e_e}\Ric,\nabla_{e_e}\Ric),\qquad \Omega=\frac{1}{\lambda}\sum_{e,i,b,d}(\nabla_{e_e}\Rm)_{ibid}(\nabla_{e_e}\Ric)_{bd},
\]
with $i$ tangent and $e,b,d$ arbitrary, and the last two terms of \eqref{eq_bochner_laplacian} collect the products of first derivatives, namely $\sum_e\big(2\Rm(\nabla_e\Ric,\nabla_e\Ric)+4(\nabla_e\Rm)(\Ric,\nabla_e\Ric)\big)$, computed with $\Ric=\lambda P_{TN}$ at $p$. The first term on the left-hand side of \eqref{eq_bochner_laplacian} is given by Lemma \ref{lem_hessian}. We now compute the three remaining curvature terms of \eqref{eq_bochner_laplacian} at $p$, and Lemmas \ref{lem_M} and \ref{lem_X} then express $\Theta$ and $\Omega$ through the geometry of $N$.

Since $\nabla f=0$ at $p$, \eqref{eq_soliton_commutation} gives $\Delta\Ric=\Delta_f\Ric=2\lambda\Ric-2\Rm(\Ric)$ at $p$. Since $\Ric=\lambda P_{TN}$, we have $\Rm(\Ric)=\lambda\tilde\rho$ and $\Rm(A,\Ric)=\lambda\sum_{a,c}\tilde\rho_{ac}A_{ac}$ for every symmetric $A$, so that, at $p$,
\begin{align}
	\Delta\Ric&=2\lambda^2P_{TN}-2\lambda\tilde\rho,\label{eq_lap_ric}\\
	2\Rm(\Delta\Ric,\Ric)&=4\lambda^3\trg\rho-4\lambda^2|\tilde\rho|^2.\label{eq_term_ricci}
\end{align}
The term $(\Delta\Rm)(\Ric,\Ric)$ contains second derivatives of $\Rm$, which the second Bianchi identity reduces to the second derivatives of $\Ric$ computed in Lemma \ref{lem_hessian}.

\begin{lemma}\label{lem_term_rm}
At $p\in N$,
\[
	(\Delta\Rm)(\Ric,\Ric)=\lambda^2\Big(2\lambda\trg\rho-2|\tilde\rho|^2+2|L^s|^2-\frac32|V|^2\Big).
\]
\end{lemma}

\begin{proof}
By the second Bianchi identity and its contraction,
\begin{align*}
	\Delta R_{abcd}&=\nabla_a(\nabla_c\Ric_{bd}-\nabla_d\Ric_{bc})-\nabla_b(\nabla_c\Ric_{ad}-\nabla_d\Ric_{ac})\\
	&\quad-\sum_e[\nabla_e,\nabla_a]R_{becd}+\sum_e[\nabla_e,\nabla_b]R_{aecd}.
\end{align*}
We contract with $\Ric\otimes\Ric=\lambda^2P_{TN}\otimes P_{TN}$, that is, we take $a=c=i$, $b=d=j$ and sum over the tangent indices $i,j$. The four second derivatives of $\Ric$ are the tangential components $\Lambda_{iijj}-\Lambda_{ijji}-\Lambda_{jiij}+\Lambda_{jjii}$, which cancel by Lemma \ref{lem_hessian}(i), since $\sigma$ is totally symmetric. Expanding the commutators by \eqref{eq_ricci_identity}, each of them produces four terms, and
\begin{align*}
	\sum_{i,j}(\Delta R)_{ijij}=\sum_{i,j,e,x}\big(&-R_{eijx}R_{xeij}-R_{eiex}R_{jxij}-R_{eiix}R_{jexj}-R_{eijx}R_{jeix}\\
	&+R_{ejix}R_{xeij}+R_{ejex}R_{ixij}+R_{ejix}R_{iexj}+R_{ejjx}R_{ieix}\big).
\end{align*}
In the second and sixth terms, $\sum_eR_{eiex}=\Ric_{ix}$, which is $\lambda\delta_{ix}$ for $x$ tangent and $0$ for $x$ normal, so that these two terms are
\[
	-\lambda\sum_{i,j}R_{jiij}=\lambda\trg\rho\qquad\text{and}\qquad\lambda\sum_{i,j}R_{ijij}=\lambda\trg\rho.
\]
In the third and eighth terms, $\sum_iR_{eiix}=-\tilde\rho_{ex}$ and $\sum_jR_{jexj}=-\tilde\rho_{ex}$, so that each of these two terms is $-|\tilde\rho|^2$. After the relabelling $i\leftrightarrow j$ in the fifth and seventh terms, the remaining four terms become
\[
	-2\sum_{i,j,e,s}R_{eijs}\big(R_{ijse}+R_{isje}\big).
\]
We split this sum according to the directions of $e$ and $s$. For fixed $s$ and tangent $e$, put $C_{eij}=R_{eijs}$. Then the sum over tangent $e,i,j$ is $\sum C_{eij}(C_{jei}-C_{ije})$, which vanishes by the cyclic relabelling $(e,i,j)\mapsto(i,j,e)$. For normal $e$ and tangent $s$, the summand is $-R_{eijs}(R_{esij}+R_{ejis})$, and the two products cancel after the exchange of $j$ and $s$. For $e=e_\alpha$ and $s=e_\beta$ normal, we have $R_{\alpha ij\beta}=-L_{\alpha\beta ij}$, $R_{ij\beta\alpha}=-V_{\alpha\beta ij}$ and $R_{i\beta j\alpha}=L_{\alpha\beta ji}$, so that this part equals $\inner{L}{V}-\inner{L}{L^t}=\frac34|V|^2-|L^s|^2$, since $\inner{L}{V}=\frac12|V|^2$, by \eqref{eq_N_norms}. Hence the four terms equal $-2\big(\frac34|V|^2-|L^s|^2\big)$. Since $(\Delta\Rm)(\Ric,\Ric)=\lambda^2\sum_{i,j}(\Delta R)_{ijij}$, the lemma follows.
\end{proof}

The term $\inner{\nabla\Ric}{\Delta\nabla\Ric}$ is computed by commuting the Laplacian with the covariant derivative.

\begin{lemma}\label{lem_term_grad}
At $p\in N$,
\[
	\inner{\nabla\Ric}{\Delta\nabla\Ric}=9\lambda^3t-2\lambda^2\Omega-6\Theta.
\]
\end{lemma}

\begin{proof}
Commuting derivatives with \eqref{eq_ricci_identity}, we obtain
\begin{align*}
	\Delta\nabla_e\Ric_{ab}&=\nabla_e\Delta\Ric_{ab}+\Ric_{ex}\nabla_x\Ric_{ab}\\
	&\quad+2\sum_{s,x}\big(R_{seax}\nabla_s\Ric_{xb}+R_{sebx}\nabla_s\Ric_{ax}\big)\\
	&\quad+\sum_x\big(\nabla_a\Ric_{ex}-\nabla_x\Ric_{ea}\big)\Ric_{xb}+\sum_x\big(\nabla_b\Ric_{ex}-\nabla_x\Ric_{eb}\big)\Ric_{ax},
\end{align*}
where $\Ric_{ex}=\sum_sR_{sesx}$ and the last line comes from the contracted second Bianchi identity $\sum_s\nabla_sR_{seax}=\nabla_a\Ric_{ex}-\nabla_x\Ric_{ea}$. The two terms in the last line vanish at $p$, since $\nabla\Ric$ is totally symmetric there by Lemma \ref{lem_first_order}(i). By \eqref{eq_soliton_commutation}, $\Delta\Ric=2\lambda\Ric-2\Rm(\Ric)+\nabla_{\nabla f}\Ric$, and, at $p$, $\nabla_e(\nabla_{\nabla f}\Ric)=\lambda(P_E)_{ex}\nabla_x\Ric$. Since $\Ric_{ex}=\lambda(P_{TN})_{ex}$, the two terms of the form $\Ric_{ex}\nabla_x\Ric_{ab}$ and $\lambda(P_E)_{ex}\nabla_x\Ric_{ab}$ combine into $\lambda\nabla_e\Ric_{ab}$, and
\begin{align*}
	\Delta\nabla_e\Ric_{ab}&=3\lambda\nabla_e\Ric_{ab}-2(\nabla_e\Rm)(\Ric)_{ab}-2\Rm(\nabla_e\Ric)_{ab}\\
	&\quad+2\sum_{s,x}\big(R_{seax}\nabla_s\Ric_{xb}+R_{sebx}\nabla_s\Ric_{ax}\big).
\end{align*}
We pair with $\nabla_e\Ric_{ab}$ and sum over $e$. By Lemma \ref{lem_first_order}(i), the first term gives $3\lambda|\nabla\Ric|^2=9\lambda^3t$. Since $\Ric=\lambda P_{TN}$, the second gives
\[
	-2\lambda\sum(\nabla_e\Rm)_{ibid}(\nabla_e\Ric)_{bd}=-2\lambda^2\Omega,
\]
and the third gives $-2\Theta$. In the fourth term, the total symmetry of $\nabla\Ric$ gives
\begin{align*}
	\sum R_{seax}\nabla_s\Ric_{xb}\nabla_e\Ric_{ab}&=\sum_bR_{seax}(\nabla\Ric)_{sxb}(\nabla\Ric)_{eab}\\
	&=-\sum_b\Rm\big((\nabla\Ric)_{\cdot\cdot b},(\nabla\Ric)_{\cdot\cdot b}\big)=-\Theta,
\end{align*}
by the skew symmetry of $R$ in its last two entries, and the last term gives $-\Theta$ in the same way.
\end{proof}

The following identity is a Bochner formula for $\nabla\Ric$ at the critical set of $f$.

\begin{proposition}\label{prop_bochner}
At every point of $N$,
\begin{align}\label{eq_bochner_point}
	\frac{1}{\lambda^2}|\nabla^2\Ric|^2+6|\rho-\lambda\Id|^2+18\lambda t+12|d|^2+\frac{27}{2}c+\frac32|V|^2-2|L^s|^2=\frac{8}{\lambda^2}\Theta+6\Omega.
\end{align}
\end{proposition}

\begin{proof}
Substituting \eqref{eq_term_ricci} and Lemmas \ref{lem_term_rm} and \ref{lem_term_grad} into \eqref{eq_bochner_laplacian}, we obtain
\[
	|\nabla^2\Ric|^2+9\lambda^3t-2\lambda^2\Omega-6\Theta=\lambda^2\Big(6\lambda\trg\rho-6|\tilde\rho|^2+2|L^s|^2-\frac32|V|^2\Big)+2\Theta+4\lambda^2\Omega.
\]
Dividing by $\lambda^2$ and using Lemma \ref{lem_hess_two}(ii) to write
\[
	6\lambda\trg\rho-6|\tilde\rho|^2-9\lambda t=-6|\rho-\lambda\Id|^2-18\lambda t-12|d|^2-\frac{27}{2}c,
\]
we obtain \eqref{eq_bochner_point}.
\end{proof}

We now express $\Theta$ and $\Omega$ through the geometry of $N$. Set
\begin{gather*}
	\omega=\sum_\alpha\sum_{i,j,\ell,m}R_{ij\ell m}h^\alpha_{i\ell}h^\alpha_{jm},\qquad A^\alpha_{\ell ij}=h^\alpha_{ij,\ell}-h^\alpha_{\ell j,i},\\
	\tau_3=-\sum_{\alpha,j,\ell}(D_{e_j}d^\alpha)_\ell h^\alpha_{j\ell},\qquad \tau_1=\sum_{\alpha,i,j,\ell}(D_{e_i}A^\alpha)_{\ell ij}h^\alpha_{j\ell},
\end{gather*}
where $d^\alpha$ and $A^\alpha$ are regarded as sections of $T^*N\otimes E$ and of $\Lambda^2T^*N\otimes T^*N\otimes E$, differentiated with the connection $D$, and $C_D=\sum_{\alpha,i,j,\ell}(D_ih^\alpha)_{j\ell}(D_jh^\alpha)_{i\ell}$.

\begin{lemma}\label{lem_M}
At $p\in N$, $\lambda^{-2}\Theta=\omega+\frac32\inner{K}{V}-\inner{P}{L^s}$.
\end{lemma}

\begin{proof}
For $e=e_\alpha$, $\nabla_\alpha\Ric$ is $\lambda S_\alpha$ on $TN$ and vanishes elsewhere, which gives $\lambda^2\omega$. For $e=e_i$, $\nabla_i\Ric$ has only the mixed entries $\lambda h^\alpha_{ij}$, and
\begin{align*}
	\Rm(\nabla_i\Ric,\nabla_i\Ric)&=\lambda^2\sum h^\alpha_{ij}h^\beta_{i\ell}\big(R_{j\ell\alpha\beta}+R_{j\beta\alpha \ell}+R_{\alpha \ell j\beta}+R_{\alpha\beta j\ell}\big)\\
	&=\lambda^2\sum h^\alpha_{ij}h^\beta_{i\ell}\big(2V_{\alpha\beta j\ell}-2L_{\alpha\beta \ell j}\big).
\end{align*}
Summing over $i$ and writing $S_\alpha S_\beta=\frac12(P_{\alpha\beta}+K_{\alpha\beta})$, we obtain
\[
	\lambda^2\sum(S_\alpha S_\beta)_{j\ell}(2V_{\alpha\beta j\ell}-2L_{\alpha\beta \ell j})=\lambda^2\big(\inner{K}{V}-\inner{P}{L^t}-\inner{K}{L^t}\big),
\]
since $P$ is symmetric and $V$ is skew in $(j,\ell)$. Finally,
\[
	\inner{P}{L^t}=\inner{P}{L^s},\qquad \inner{K}{L^t}=\inner{K}{L^s-\tfrac12V}=-\tfrac12\inner{K}{V}.
\]
\end{proof}

\begin{lemma}\label{lem_X}
At $p\in N$,
\[
	\Omega=3\tau_3-\tau_1+3\inner{T}{\rho}+\omega+2\inner{P}{L^s}-\frac32\inner{K}{V}+\frac92c.
\]
On the closed manifold $N$,
\begin{align}\label{eq_tau_integrals}
	\int_N\tau_3\dmu_N=\int_N|d|^2\dmu_N,\qquad \int_N\tau_1\dmu_N=\int_N\big(|D\II|^2-C_D\big)\dmu_N.
\end{align}
\end{lemma}

\begin{proof}
By the definition of $\Omega$ and Lemma \ref{lem_first_order}(i),
\[
	\Omega=\sum(\nabla_\alpha\Rm)_{iji\ell}h^\alpha_{j\ell}+2\sum(\nabla_\ell\Rm)_{i\alpha ij}h^\alpha_{\ell j}.
\]
The second Bianchi identity $\nabla_\alpha R_{iji\ell}=-\nabla_iR_{j\alpha i\ell}+\nabla_jR_{i\alpha i\ell}$ applied to the first sum gives $\Omega=3\Omega_3-\Omega_1$, with
\[
	\Omega_3=\sum(\nabla_j\Rm)_{i\alpha i\ell}h^\alpha_{j\ell},\qquad \Omega_1=\sum(\nabla_i\Rm)_{j\alpha i\ell}h^\alpha_{j\ell}.
\]
We expand the covariant derivatives in the frames of the proof of Lemma \ref{lem_hessian}. For $\Omega_3$, the components summed over $i$ are $-d^\alpha_\ell$ by Lemma \ref{lem_first_order}(ii), whose derivative gives $\tau_3$, and for $\Omega_1$, the components are $R_{j\alpha i\ell}=A^\alpha_{\ell ij}$ by Lemma \ref{lem_first_order}(ii), whose derivative gives $\tau_1$. The eight contractions of the connection terms are computed in Appendix \ref{app_conn}, and give
\begin{align}\label{eq_Omega}
	\Omega_3=\tau_3+\inner{T}{\rho}+\frac12\inner{P}{L^s}-\frac34\inner{K}{V}+\frac32c,\qquad \Omega_1=\tau_1-\omega-\frac12\inner{P}{L^s}-\frac34\inner{K}{V},
\end{align}
which gives the formula for $\Omega$. For \eqref{eq_tau_integrals}, consider the vector fields on $N$
\[
	Z_3=\sum_{\alpha,j,\ell}d^\alpha_\ell h^\alpha_{j\ell}\,e_j,\qquad Z_1=\sum_{\alpha,i,j,\ell}A^\alpha_{\ell ij}h^\alpha_{j\ell}\,e_i,
\]
which do not depend on the choice of the orthonormal frames. Since $\sum_j h^\alpha_{j\ell,j}=d^\alpha_\ell$, their divergences at $p$ are
\[
	\divg Z_3=\sum(D_{e_j}d^\alpha)_\ell h^\alpha_{j\ell}+\sum d^\alpha_\ell d^\alpha_\ell=-\tau_3+|d|^2
\]
and
\[
	\divg Z_1=\sum(D_{e_i}A^\alpha)_{\ell ij}h^\alpha_{j\ell}+\sum A^\alpha_{\ell ij}h^\alpha_{j\ell,i}=\tau_1+C_D-|D\II|^2,
\]
where, in the last equality, we used $A^\alpha_{\ell ij}=h^\alpha_{ij,\ell}-h^\alpha_{\ell j,i}$ and
\[
	\sum h^\alpha_{ij,\ell}h^\alpha_{j\ell,i}=C_D,\qquad \sum h^\alpha_{\ell j,i}h^\alpha_{j\ell,i}=|D\II|^2.
\] Integrating over the closed manifold $N$ gives \eqref{eq_tau_integrals}.
\end{proof}

The tensor $A^\alpha$ is skew in $(\ell,i)$, and $|D\II|^2-C_D=\frac12|A|^2$, since $|A|^2=2|D\II|^2-2C_D$. Moreover,
\begin{align}\label{eq_d_bound}
	|d|^2\leq\frac{k-1}{2}|A|^2,
\end{align}
with equality when $k=2$: since $\II$ is trace free, $d^\alpha_j=\sum_iA^\alpha_{iji}$, and the part of $A^\alpha$ of the form $\delta_{ij}u_\ell-\delta_{\ell j}u_i$ is orthogonal to the trace free part and has norm $\frac{2}{k-1}|d^\alpha|^2$; when $k=2$ the space of tensors skew in $(\ell,i)$ has dimension two per normal direction, and the trace part is everything.

The identity of Theorem \ref{thm_general} also contains the quantity $C=\frac12\sum_{\alpha,\beta}|K_{\alpha\beta}|^2$, and we write $a=\sum_{\alpha,\beta}\trg(S_\alpha S_\beta S_\alpha S_\beta)$. Then one checks directly that
\begin{align}\label{eq_algebraic_norms}
	a=|T|^2-C,\qquad |P|^2=4|T|^2-2C,\qquad |K|^2=2C,\qquad |\sigma|^2=3c+6a.
\end{align}

\begin{theorem}\label{thm_general}
For every $k\geq2$ and $q\geq1$, $\int_NJ\dmu_N=0$, where
\begin{align*}
	J={}&4|D\II|^2+3|A|^2-6|d|^2+4\Big|L^s+\frac12P\Big|^2+2\Big|V-\frac34K\Big|^2+2|Y|^2+\frac{1}{\lambda^2}|F|^2\\
	&+6|\rho-\lambda\Id|^2+18\lambda t+|\sigma|^2+5|P|^2-\frac94C-\frac{27}{2}c-14\omega-18\inner{T}{\rho}.
\end{align*}
\end{theorem}

\begin{proof}
Insert \eqref{eq_hessian_norm} and Lemmas \ref{lem_M} and \ref{lem_X} into \eqref{eq_bochner_point}. The result is the pointwise identity
\begin{align}\label{eq_pointwise}
\begin{split}
	0={}&|\sigma|^2+4|D\II|^2+12|d|^2+6|P|^2+4\inner{P}{L^s}+4|L^s|^2+2|V|^2-3\inner{K}{V}+2|Y|^2+\frac{1}{\lambda^2}|F|^2\\
	&+6|\rho-\lambda\Id|^2+18\lambda t-\frac{27}{2}c-14\omega-18\inner{T}{\rho}-18\tau_3+6\tau_1.
\end{split}
\end{align}
Since $6|P|^2+4\inner{P}{L^s}+4|L^s|^2=4|L^s+\frac12P|^2+5|P|^2$ and
\[
	2|V|^2-3\inner{K}{V}=2\Big|V-\frac34K\Big|^2-\frac98|K|^2=2\Big|V-\frac34K\Big|^2-\frac94C,
\]
we integrate over $N$ and use \eqref{eq_tau_integrals} to replace $-18\tau_3+6\tau_1$ by $-18|d|^2+6|D\II|^2-6C_D$. Finally, $6(|D\II|^2-C_D)=3|A|^2$.
\end{proof}

The last line of $J$ can be expressed through the intrinsic geometry of $N$, which gives the form of $J$ used in Sections \ref{sec_lowdim} and \ref{sec_four}. By Lemma \ref{lem_hess_two}(ii), the trace free part of $\Ric_N$ is $\mathring{\Ric}_N=\Ric_N-(\lambda+\frac{t}{2k})g_N$, and $\rho-\lambda\Id=\frac{3t}{2k}\Id+\mathring{\Ric}_N+\mathring{T}$, where $\mathring{T}=T-\frac{t}{k}\Id$ is the trace free part of $T$. The Gauss equation gives $\omega=\sum_\alpha R^N(S_\alpha,S_\alpha)-c+a$, where $R^N(S,S)=\sum R^N_{ij\ell m}S_{i\ell}S_{jm}$, and, when $k\geq3$, the decomposition of $R^N$ into its Weyl and Schouten parts gives $R^N(S,S)=W^N(S,S)-2\inner{\mathrm{Sch}^N}{S^2}$ for trace free $S$, where
\[
	\mathrm{Sch}^N=\frac{1}{k-2}\Big(\Ric_N-\frac{R_N}{2(k-1)}g_N\Big)=\frac{1}{k-2}\mathring{\Ric}_N+\frac{\lambda}{2(k-1)}g_N+\frac{t}{4k(k-1)}g_N,
\]
so that
\begin{align}\label{eq_schouten_pairing}
	\sum_\alpha R^N(S_\alpha,S_\alpha)=W^N(\II,\II)-\frac{2}{k-2}\inner{\mathring{\Ric}_N}{\mathring{T}}-\frac{\lambda t}{k-1}-\frac{t^2}{2k(k-1)},
\end{align}
where $W^N(\II,\II)=\sum_\alpha W^N(S_\alpha,S_\alpha)$.
Substituting \eqref{eq_algebraic_norms} and \eqref{eq_schouten_pairing} into $J$, the terms $|\mathring{T}|^2$ cancel, and we obtain, for $k\geq3$,
\begin{align}\label{eq_general_weyl}
\begin{split}
	J={}&4|D\II|^2+3|A|^2-6|d|^2+4\Big|L^s+\frac12P\Big|^2+2\Big|V-\frac34K\Big|^2+2|Y|^2+\frac{1}{\lambda^2}|F|^2\\
	&+\frac{14}{k-1}\lambda t+6|\mathring{\Ric}_N|^2+\Big(\frac{28}{k-2}-6\Big)\inner{\mathring{\Ric}_N}{\mathring{T}}-14W^N(\II,\II)\\
	&-\frac{17}{4}C+\frac72c+\Big(\frac{7}{k(k-1)}-\frac{3}{2k}\Big)t^2.
\end{split}
\end{align}
When $k=2$, $\Ric_N$ and $S_\alpha^2$ are multiples of the identity, so that $\mathring{\Ric}_N=0$ and $\mathring{T}=0$, and $R^N(S,S)=-K_N|S|^2$ for trace free $S$, with $K_N=\lambda+\frac14t$ by Lemma \ref{lem_hess_two}, so that $\sum_\alpha R^N(S_\alpha,S_\alpha)=-\lambda t-\frac14t^2$; hence \eqref{eq_general_weyl} holds with the terms containing $\mathring{\Ric}_N$, $\mathring{T}$ and $W^N$ omitted.

We also need the integration by parts identity for the second fundamental form which follows from the Ricci identity for $D$. Multiples of it are subtracted from $J$ in the proofs of Theorems \ref{thm_three} and \ref{thm_four}. Set
\begin{gather*}
	O_{\II}=\sum_{i,j,\ell}\inner{R^\perp(e_i,e_j)\II_{j\ell}}{\II_{i\ell}}=\frac12\inner{K}{V}+C,\\
	\omega_{\II}=\sum_{i,j,\ell,m}R^N_{ij\ell m}\inner{\II_{im}}{\II_{j\ell}}=-\sum_\alpha R^N(S_\alpha,S_\alpha),
\end{gather*}
where the value of $O_{\II}$ follows from Lemma \ref{lem_first_order}(iii), since $\sum_\ell h^\alpha_{i\ell}h^\beta_{j\ell}=\frac12(P_{\alpha\beta}+K_{\alpha\beta})_{ij}$.

\begin{lemma}\label{lem_simons}
On the closed manifold $N$, $\int_NI\dmu_N=0$, where
\[
	I=C_D-|d|^2-O_{\II}+\trg(\Ric_NT)+\omega_{\II}.
\]
For every $k\geq3$,
\begin{align}\label{eq_simons_weyl}
	\trg(\Ric_NT)+\omega_{\II}=-W^N(\II,\II)+\frac{k}{k-2}\inner{\mathring{\Ric}_N}{\mathring{T}}+\frac{k}{k-1}\lambda t+\frac{t^2}{2(k-1)}.
\end{align}
\end{lemma}

\begin{proof}
Regard $\II$ as the one-form $\beta\in\Omega^1(N;T^*N\otimes E)$, $\beta(X)(Y)=\II(X,Y)$, with values in the bundle $T^*N\otimes E$ with the connection $D$. Then $d^D\beta$ has components $A^\alpha_{\ell ij}$, so that $|d^D\beta|^2=\frac12|A|^2$ in the norm of two-forms, and $\delta^D\beta=-d$ by the definition of $d$. The integrated Weitzenb\"ock formula for bundle-valued one-forms \cite[Proposition A.7]{drth} gives
\[
	\int_N\Big(\frac12|A|^2+|d|^2\Big)\dmu_N=\int_N\big(|D\II|^2-O_{\II}+\trg(\Ric_NT)+\omega_{\II}\big)\dmu_N.
\]
Here the curvature term of the formula is $\sum_{i,j,\ell}\inner{\II_{j\ell}}{\big((D^2_{e_i,e_j}-D^2_{e_j,e_i})\II\big)_{i\ell}}$, and, by the Ricci identity for $D$, the curvature $R^\perp$ of the normal connection contributes $-O_{\II}$ to it, while the curvature $R^N(X,Y)=\nabla^N_X\nabla^N_Y-\nabla^N_Y\nabla^N_X-\nabla^N_{[X,Y]}$, acting on the first and on the second entry of $\II$, contributes $\trg(\Ric_NT)$ and $\omega_{\II}$. Since $\frac12|A|^2=|D\II|^2-C_D$, this is $\int_NI\dmu_N=0$. Finally,
\[
	\trg(\Ric_NT)=\lambda t+\frac{t^2}{2k}+\inner{\mathring{\Ric}_N}{\mathring{T}},
\]
and \eqref{eq_schouten_pairing} gives \eqref{eq_simons_weyl}.
\end{proof}

\section{Two- and three-dimensional focal sets}\label{sec_lowdim}

We now prove Theorem \ref{thm_focal} and Corollary \ref{cor_alldim}. When $k=2$, the identity of Theorem \ref{thm_general} suffices, and, when $k=3$, we combine it with Lemma \ref{lem_simons}.

\begin{theorem}\label{thm_two}
If $k=2$, then, for every $q\geq1$,
\begin{align}\label{eq_two}
\begin{split}
	0=\int_N\Big(&4|D\II|^2+14\lambda t+\frac{19}{8}t^2+\frac{31}{4}\Big(c-\frac{t^2}{2}\Big)\\
	&+4\Big|L^s+\frac12P\Big|^2+2\Big|V-\frac34K\Big|^2+2|Y|^2+\frac{1}{\lambda^2}|F|^2\Big)\dmu_N,
\end{split}
\end{align}
and every term of the integrand is nonnegative. Hence $N$ is totally geodesic.
\end{theorem}

\begin{proof}
For trace free symmetric endomorphisms of a two-dimensional space, polarization of $S^2=\frac12|S|^2\Id$ gives $S_\alpha S_\beta+S_\beta S_\alpha=G_{\alpha\beta}\Id$. Hence
\[
	T=\frac t2\Id,\qquad P_{\alpha\beta}=G_{\alpha\beta}\Id,\qquad |P|^2=2c.
\]
Since $|T|^2=\frac12t^2$ and $|P|^2=2c$, \eqref{eq_algebraic_norms} gives $C=t^2-c$. Moreover, $3|A|^2-6|d|^2=0$ by the equality case of \eqref{eq_d_bound}. Inserting these into the form of \eqref{eq_general_weyl} for $k=2$, we obtain
\begin{align*}
	J&=4|D\II|^2+(\text{the four squares})+14\lambda t+\frac{31}{4}c-\frac32t^2\\
	&=4|D\II|^2+(\text{the four squares})+14\lambda t+\frac{19}{8}t^2+\frac{31}{4}\Big(c-\frac{t^2}{2}\Big),
\end{align*}
which is \eqref{eq_two}. The matrix $G$ is the Gram matrix of $S_1,\ldots,S_q$ in the two-dimensional Euclidean space of trace free symmetric endomorphisms of $T_pN$, so $G\geq0$ has rank at most two and trace $t$, whence $c=|G|^2\geq\frac12t^2$. Every term of the integrand of \eqref{eq_two} is therefore nonnegative, and since $\lambda>0$, the term $14\lambda t$ forces $t=0$, so $N$ is totally geodesic.
\end{proof}

\begin{theorem}\label{thm_three}
If $k=3$, then, for every $q\geq1$,
\begin{align}\label{eq_three}
\begin{split}
	0=\int_N\Big(&5|A|^2-2|d|^2+\lambda t+6\Big|\mathring{\Ric}_N+\frac56\mathring{T}\Big|^2+\frac{13}{36}\big(C+t^2\big)+\frac{19}{9}c\\
	&+4\Big|L^s+\frac12P\Big|^2+2\Big|V-\frac14K\Big|^2+2|Y|^2+\frac{1}{\lambda^2}|F|^2\Big)\dmu_N,
\end{split}
\end{align}
and every term of the integrand is nonnegative, since $|d|^2\leq|A|^2$ by \eqref{eq_d_bound}. Hence $N$ is totally geodesic.
\end{theorem}

\begin{proof}
Subtract four times the vanishing integral $\int_NI\dmu_N$ of Lemma \ref{lem_simons} from the identity of Theorem \ref{thm_general}, written in the form \eqref{eq_general_weyl}. The derivative terms become
\[
	4|D\II|^2+3|A|^2-6|d|^2-4(C_D-|d|^2)=5|A|^2-2|d|^2,
\]
since $|D\II|^2-C_D=\frac12|A|^2$. The term $4O_{\II}=2\inner{K}{V}+4C$ combines with $2|V-\frac34K|^2-\frac{17}{4}C$ into $2|V-\frac14K|^2+\frac74C$. By \eqref{eq_simons_weyl} with $k=3$, where $W^N=0$, the remaining terms become
\[
	\lambda t+6|\mathring{\Ric}_N|^2+10\inner{\mathring{\Ric}_N}{\mathring{T}}+\frac74C+\frac72c-\frac13t^2.
\]
Finally,
\[
	6|\mathring{\Ric}_N|^2+10\inner{\mathring{\Ric}_N}{\mathring{T}}=6\Big|\mathring{\Ric}_N+\frac56\mathring{T}\Big|^2-\frac{25}{6}|\mathring{T}|^2,
\]
and the Cayley--Hamilton identity $\trg S^4=\frac12(\trg S^2)^2$ for trace free symmetric $3\times3$ matrices, applied to $S=\sum_\alpha z^\alpha S_\alpha$ and polarized in $z$, gives
\[
	2|T|^2+a=c+\frac{t^2}{2},\qquad\text{that is,}\qquad |\mathring{T}|^2=\frac13\Big(C+c-\frac{t^2}{2}\Big),
\]
by \eqref{eq_algebraic_norms}. Substituting gives \eqref{eq_three}, and the term $\lambda t$ forces $t=0$.
\end{proof}

\begin{proof}[Proof of Theorem \ref{thm_focal}]
By Theorems \ref{thm_two} and \ref{thm_three}.
\end{proof}

\begin{proof}[Proof of Corollary \ref{cor_alldim}]
If $f$ is constant, then $n=k\in\{2,3\}$ and $\Ric=\lambda g$, so $M$ is a quotient of $\Sph^k(\sqrt{(k-1)/\lambda})$ by a finite group. Otherwise, $N$ has dimension $k\in\{2,3\}$ by Proposition \ref{prop_structure}, and it is totally geodesic by Theorem \ref{thm_focal}. By \cite[Corollary 1.1]{lloqwz}, $M$ is rigid, so its universal cover is $\widetilde N\times\R^{n-k}$, where $\widetilde N$ is a complete simply connected Einstein manifold of dimension $k$ with $\Ric_{\widetilde N}=\lambda g_{\widetilde N}$, that is, the sphere $\Sph^k(\sqrt{(k-1)/\lambda})$, the potential function is $\frac\lambda2|x|^2$ on the Euclidean factor up to a constant, and the deck group $\Gamma$ acts freely on $\widetilde N$ and orthogonally on $\R^{n-k}$. Since $\Gamma$ acts properly discontinuously and preserves the compact set $\widetilde N\times\{0\}$, it is finite. Finally, $R=0$ gives the Gaussian soliton and $R=n\lambda$ a compact Einstein manifold \cite[Propositions 3.1 and 3.3]{pewy}, while $R=\lambda$ does not occur and $R=(n-1)\lambda$ gives rigid solitons \cite[Theorem 10]{fr16}, which leaves $R\in\{2\lambda,3\lambda\}$ when $R\leq3\lambda$, and $R=4\lambda$ when $n=6$.
\end{proof}

\begin{remark}[K\"ahler shrinkers]\label{rem_kahler}
Suppose that $(M,g)$ is K\"ahler with complex structure $J$ and satisfies \eqref{fundeq}. Then $\Ric$ and $\nabla^2f$ are $J$-invariant, so $\nabla f$ is real holomorphic and $J\nabla f$ is a Killing field. Its zero set is the critical set of $f$, which is $N$ by Proposition \ref{prop_structure}. By a theorem of Kobayashi \cite{kobayashi}, each connected component of the zero set of a Killing field is a totally geodesic submanifold of even codimension. Hence $N$ is totally geodesic for every complete K\"ahler gradient shrinking Ricci soliton with constant scalar curvature, and \cite[Corollary 1.1]{lloqwz} gives a further route to the rigidity of these solitons, established in \cite{hloqwz,horkahler,lowz}.
\end{remark}

\section{Totally geodesic minimum sets in dimension five}\label{sec_level}

In this section and the next, $(M^n,g,f,\lambda)$ is a complete gradient shrinking Ricci soliton with constant scalar curvature $R=k\lambda$ and nonconstant potential function, and we use the following standard facts on the level sets of $f$. On $M\setminus N$, let $E_1=\nabla r$, let $\Sigma_r=r^{-1}(r)$ be the level set of $r$ through a given point, and let $SX=\nabla_XE_1$ be its shape operator, with mean curvature $H=\trg S$. Since $|\nabla r|=1$, we have $\nabla_{E_1}E_1=0$, so that the integral curves of $E_1$ are the normal geodesics of $N$, and we use parallel orthonormal frames along them, tangent to the level sets, to identify tensors on different level sets. For a tensor $T$ on the level sets, we write $T'=\nabla_{E_1}T$, which is the derivative of the matrix of $T$ in such a frame. We denote by $\Ric_T$ the restriction of $\Ric$ to $E_1^\perp=T\Sigma_r$, and by $\nabla^\Sigma$, $\Rm_\Sigma$, $\Ric_\Sigma$, $R_\Sigma$ and $W_\Sigma$ the Levi-Civita connection, the curvature tensor, the Ricci operator, the scalar curvature and the Weyl tensor of the metric induced on $\Sigma_r$.

By Proposition \ref{prop_structure}, every level set $\Sigma_r$, $r>0$, is the image under $\exp$ of the bundle of normal vectors of length $r$, so it is a compact hypersurface, the total space of a sphere bundle over $N$ with fiber $\Sph^{q-1}$, which is connected when $q\geq2$, and $M$ is homotopy equivalent to $N$. Since the volume measure of $M$ is $\mathrm{d}\mu_N\,\mathrm{d}z$ and $|\nabla r|=1$, the volume measure of $\Sigma_r$ is $r^{q-1}\,\mathrm{d}\sigma$ under the identification with the unit normal bundle of $N$, where $\mathrm{d}\sigma$ is the product of the measure of $N$ and of the standard measure of the unit spheres of $E$. Since $\Ric(E_1,\cdot)=0$ and $\nabla f=\lambda rE_1$, equation \eqref{fundeq} gives $\nabla^2f(X,Y)=\lambda r\inner{SX}{Y}$ for $X,Y$ tangent to $\Sigma_r$, that is,
\begin{align}\label{eq_ricT}
	\Ric_T=\lambda\Id-\lambda rS,\qquad \trg\Ric_T=k\lambda,\qquad |\Ric_T|^2=k\lambda^2,
\end{align}
where the last two identities follow from $\trg\Ric=R$ and $|\Ric|^2=\lambda R$. Equivalently,
\begin{align}\label{eq_rH}
	rH=r^2|S|^2=q-1.
\end{align}
The Riccati equation for the unit geodesic field $E_1$ is
\[
	\Rm(E_1,Y,E_1,T)=-\inner{(S'+S^2)Y}{T},\qquad Y,T\perp E_1,
\]
and, taking traces in the Gauss equation, we obtain $\Ric_\Sigma=\Ric_T-\Rm(E_1,\cdot,E_1,\cdot)+HS-S^2$, that is,
\begin{align}\label{eq_ric_sigma}
	\Ric_T=\Ric_\Sigma-HS-S'.
\end{align}
Comparing with \eqref{eq_ricT}, we obtain
\begin{align}\label{eq_Sprime}
	S'=\Ric_\Sigma-\lambda\Id+(\lambda r-H)S,
\end{align}
see also \cite[Lemma 3.1]{horlevel}. We first restrict the Bochner formula \eqref{eq_bochner_global} to the level sets. We write $\Ric_T'=(\nabla_{E_1}\Ric)_T$ for the restriction of $\nabla_{E_1}\Ric$ to $T\Sigma_r$.

\begin{lemma}\label{lem_bochner_level}
On every level set $\Sigma_r$, $r>0$, we have $\divg_\Sigma\Ric_T=0$ and
\begin{align}\label{eq_bochner_level}
	|\nabla^\Sigma\Ric_T|^2+|\Ric_T'|^2=2\Rm_\Sigma(\Ric_T,\Ric_T)-2k\lambda^3.
\end{align}
\end{lemma}

\begin{proof}
Let $X,Y,Z$ be tangent to $\Sigma_r$. Since $\Ric(E_1,\cdot)=0$, $\nabla_{E_1}E_1=0$ and $\nabla_XY=\nabla^\Sigma_XY-\inner{SX}{Y}E_1$, we have
\begin{gather*}
	(\nabla_X\Ric)(Y,Z)=(\nabla^\Sigma_X\Ric_T)(Y,Z),\qquad (\nabla_{E_1}\Ric)(Y,Z)=\Ric_T'(Y,Z),\\
	(\nabla_X\Ric)(E_1,Y)=-\Ric(SX,Y),
\end{gather*}
while $(\nabla_{E_1}\Ric)(E_1,\cdot)$ and $(\nabla_X\Ric)(E_1,E_1)$ vanish. Hence
\begin{align}\label{eq_gradric_split}
	|\nabla\Ric|^2=|\nabla^\Sigma\Ric_T|^2+|\Ric_T'|^2+2|\Ric_TS|^2.
\end{align}
Moreover, if $\{e_a\}$ is an orthonormal frame of $T\Sigma_r$, then
\begin{align*}
	(\divg_\Sigma\Ric_T)(Y)&=\sum_a(\nabla_{e_a}\Ric)(e_a,Y)\\
	&=(\divg\Ric)(Y)-(\nabla_{E_1}\Ric)(E_1,Y)=\frac12\nabla_YR=0.
\end{align*}

By \eqref{eq_ricT}, $S$ commutes with $\Ric_T$, and
\[
	\trg(\Ric_TS)=\frac{1}{\lambda r}\big(\lambda\trg\Ric_T-|\Ric_T|^2\big)=0.
\]
Let $\{e_a\}$ be an orthonormal basis of $T_x\Sigma_r$ of eigenvectors of $\Ric_T$ and $S$, with $\Ric_Te_a=\rho_ae_a$ and $Se_a=\kappa_ae_a$. By the Gauss equation, $\Rm(e_a,e_b,e_a,e_b)=\Rm_\Sigma(e_a,e_b,e_a,e_b)-\kappa_a\kappa_b$ for $a\neq b$, and, since $\Ric(E_1,\cdot)=0$,
\begin{align*}
	\Rm(\Ric,\Ric)&=\Rm_\Sigma(\Ric_T,\Ric_T)-\sum_{a\neq b}\kappa_a\kappa_b\rho_a\rho_b\\
	&=\Rm_\Sigma(\Ric_T,\Ric_T)-\trg(\Ric_TS)^2+|\Ric_TS|^2.
\end{align*}
Since $2\lambda|\Ric|^2=2k\lambda^3$, substituting this and \eqref{eq_gradric_split} into \eqref{eq_bochner_global} gives \eqref{eq_bochner_level}.
\end{proof}

For the remainder of this section, $n=5$ and $R=2\lambda$, so that $k=2$ and $q=3$, and we prove Theorems \ref{thm_identity} and \ref{thm_level}. The level sets $\Sigma_r$ are four-dimensional, and they are $2$-sphere bundles over $N$, so that
\begin{align}\label{eq_level_volume}
	|\Sigma_r|=4\pi r^2|N|,\qquad \chi(\Sigma_r)=2\chi(N),
\end{align}
where the second identity is the multiplicativity of the Euler characteristic in the sphere bundle $\Sigma_r\to N$. By \eqref{eq_gauss_bonnet_N}, $\chi(N)>0$, and $N$ is diffeomorphic to $\Sph^2$ or to the real projective plane. By \eqref{eq_ricT} and \eqref{eq_rH},
\begin{align}\label{eq_ricT_five}
	\Ric_T=\lambda\Id-\lambda rS,\qquad \trg\Ric_T=2\lambda,\qquad |\Ric_T|^2=2\lambda^2,\qquad rH=r^2|S|^2=2.
\end{align}

\begin{lemma}\label{lem_ric_sigma_five}
On $M\setminus N$, we have
\begin{align}\label{eq_ric_sigma_five}
	\Ric_\Sigma=\frac{1}{r^2}\Id+\Big(1-\frac{1}{\lambda r^2}\Big)\Ric_T-\frac{1}{\lambda r}\Ric_T'.
\end{align}
Consequently, with $\mathring{\Ric}_T=\Ric_T-\frac\lambda2\Id$,
\begin{align}\label{eq_traceless_sigma}
\begin{split}
	R_\Sigma&=2\lambda+\frac{2}{r^2},\qquad \mathring{\Ric}_\Sigma=\beta_r\mathring{\Ric}_T-\gamma_r\Ric_T',\\
	|\mathring{\Ric}_\Sigma|^2&=\Big(\lambda-\frac{1}{r^2}\Big)^2+\frac{1}{\lambda^2r^2}|\Ric_T'|^2,
\end{split}
\end{align}
where $\mathring{\Ric}_\Sigma=\Ric_\Sigma-\frac{R_\Sigma}{4}\Id$ is the trace free Ricci tensor of $\Sigma_r$, $\beta_r=1-\frac{1}{\lambda r^2}$ and $\gamma_r=\frac{1}{\lambda r}$. Moreover,
\[
	\trg\mathring{\Ric}_T=0,\qquad |\mathring{\Ric}_T|^2=\lambda^2,\qquad \trg\Ric_T'=\inner{\Ric_T'}{\mathring{\Ric}_T}=0,\qquad \divg_\Sigma\mathring{\Ric}_T=0.
\]
\end{lemma}

\begin{proof}
Substituting $S=(\lambda r)^{-1}(\lambda\Id-\Ric_T)$ into \eqref{eq_Sprime}, and using $\lambda r-H=\lambda r-2r^{-1}$ from \eqref{eq_ricT_five}, gives \eqref{eq_ric_sigma_five}. Tracing \eqref{eq_ric_sigma_five} and using $\trg\Ric_T'=\inner{\Ric_T}{\Ric_T'}=0$, which follow from \eqref{eq_ricT_five}, gives the remaining identities. Finally, $\divg_\Sigma\mathring{\Ric}_T=\divg_\Sigma\Ric_T=0$ by Lemma \ref{lem_bochner_level}.
\end{proof}

With $\lambda=\frac12$, the formulas of Lemma \ref{lem_ric_sigma_five} for $R_\Sigma$, $\Ric_\Sigma$ and the norm of the Ricci tensor of $\Sigma_r$ also appear, in an equivalent form, in \cite[(5.16)--(5.18)]{loqw2}.

We now recall some facts about four-dimensional curvature. Let $(\Sigma^4,h)$ be a Riemannian four-manifold. We regard an algebraic curvature tensor $T$ as the self-adjoint endomorphism of $\Lambda^2T\Sigma$ given by $\inner{T(e_i\wedge e_j)}{e_m\wedge e_\ell}=T(e_i,e_j,e_m,e_\ell)$, where $\Lambda^2T\Sigma$ carries the inner product for which $\{e_i\wedge e_j\}_{i<j}$ is orthonormal when $\{e_i\}$ is orthonormal. With this convention, the curvature operator of the round sphere of radius one is the identity, and we write $\|T\|$ for the norm of an algebraic curvature tensor $T$ as an endomorphism of $\Lambda^2T\Sigma$. For the Weyl tensor $W$,
\[
	\|W\|^2=\frac14\sum_{i,j,m,\ell}W(e_i,e_j,e_m,e_\ell)^2=\frac14|W|^2,
\]
where $|W|$ is the full tensor norm. If $\Sigma$ is closed, the Chern--Gauss--Bonnet formula reads
\begin{align}\label{eq_cgb}
	8\pi^2\chi(\Sigma)=\int_\Sigma\Big(\|W\|^2-\frac12|\mathring{\Ric}|^2+\frac{R^2}{24}\Big)\dmu,
\end{align}
see \cite[Chapter 6]{besse}. It also holds when $\Sigma$ is not orientable, as follows by passing to the orientation covering. Since
\[
	\Rm(e_i,e_j,e_i,e_j)=W(e_i,e_j,e_i,e_j)+\frac12\big(\mathring{\Ric}(e_i,e_i)+\mathring{\Ric}(e_j,e_j)\big)+\frac{R}{12}
\]
for orthonormal $e_i,e_j$, a computation in an orthonormal basis of eigenvectors of a self-adjoint endomorphism $A$ gives
\begin{align}\label{eq_four_decomp}
	\Rm(A,A)=W(A,A)+\trg\big(\mathring{\Ric}A((\trg A)\Id-A)\big)+\frac{R}{12}\big((\trg A)^2-|A|^2\big),
\end{align}
where
\[
	T(A,A)=\sum_{i,j,m,\ell}T(e_i,e_j,e_m,e_\ell)\inner{Ae_i}{e_m}\inner{Ae_j}{e_\ell}.
\]
For a trace free self-adjoint endomorphism $B$ of $T_x\Sigma$, let $B\wedge B$ be the algebraic curvature tensor $(B\wedge B)(X,Y,Z,T)=\inner{BX}{Z}\inner{BY}{T}-\inner{BX}{T}\inner{BY}{Z}$, and let $Z_B=2\,\Pi_W(B\wedge B)$ be twice its Weyl part. If $\{e_i\}$ is an orthonormal basis of eigenvectors of $B$ with eigenvalues $\beta_i$, then $B\wedge B$, its Ricci contraction $-B^2$ and its scalar curvature $-|B|^2$ are diagonal in this basis, so $Z_B$ has only sectional components, and these are
\begin{align}\label{eq_ZB}
	Z_B(e_i,e_j,e_i,e_j)=(\beta_i+\beta_j)^2-\frac13|B|^2,\qquad i\neq j.
\end{align}
Since $\sum_i\beta_i=0$, the sums $\beta_i+\beta_j$ over the three pairs $\{i,j\}$ containing a fixed index have squared norm $|B|^2$, and a direct computation gives
\begin{align}\label{eq_ZB_norms}
	\inner{W}{Z_B}_{\mathrm{HS}}=W(B,B),\qquad \|Z_B\|^2=\sum_{i<j}Z_B(e_i,e_j,e_i,e_j)^2=\frac43|B|^4-4\Big|B^2-\frac{|B|^2}{4}\Id\Big|^2,
\end{align}
for every algebraic Weyl tensor $W$ of $T_x\Sigma$.

\begin{proposition}\label{prop_weyl_norm}
For every $r>0$,
\begin{align}\label{eq_weyl_norm}
	\int_{\Sigma_r}\Big(\|W_\Sigma\|^2-\frac{R_\Sigma^2}{12}\Big)\dmu_r=2\pi\int_N|\II|^2\dmu_N+\frac{1}{2\lambda^2r^2}\int_{\Sigma_r}|\Ric_T'|^2\dmu_r.
\end{align}
\end{proposition}

\begin{proof}
By \eqref{eq_level_volume}, \eqref{eq_traceless_sigma} and the Chern--Gauss--Bonnet formula \eqref{eq_cgb}, since $R_\Sigma$ is constant on $\Sigma_r$,
\[
	\int_{\Sigma_r}\Big(\|W_\Sigma\|^2-\frac{R_\Sigma^2}{12}\Big)\dmu_r=16\pi^2\chi(N)-8\pi\lambda|N|+\frac{1}{2\lambda^2r^2}\int_{\Sigma_r}|\Ric_T'|^2\dmu_r.
\]
By \eqref{eq_gauss_bonnet_N}, $16\pi^2\chi(N)-8\pi\lambda|N|=2\pi\int_N|\II|^2\dmu_N$, which gives \eqref{eq_weyl_norm}.
\end{proof}

We also need the following integral identity for divergence free symmetric $2$-tensors, which we apply to $\mathring{\Ric}_T$ on $\Sigma_r$ in the proof of Theorem \ref{thm_identity}.

\begin{lemma}\label{lem_codazzi_int}
Let $(\Sigma,h)$ be a closed Riemannian manifold and let $B$ be a symmetric $2$-tensor on $\Sigma$ with $\divg B=0$. Let $\Phi_B(X,Y,Z)=(\nabla_XB)(Y,Z)-(\nabla_YB)(X,Z)$. Then
\[
	\int_\Sigma|\Phi_B|^2\dmu=2\int_\Sigma\big(|\nabla B|^2+\trg(\Ric B^2)-\Rm(B,B)\big)\dmu.
\]
\end{lemma}

\begin{proof}
Apply the integrated first Weitzenb\"ock formula $(\delta^\nabla d^\nabla+d^\nabla\delta^\nabla)B=\nabla^*\nabla B-\mathring{R}B+B\circ\Ric$ \cite[12.69]{besse} (see also \cite[Proposition 4.1]{bourguignon}) to $B$, regarded as a $T^*\Sigma$-valued one-form, using $\delta^\nabla B=0$, $\inner{\mathring{R}B}{B}=\Rm(B,B)$ and $|\Phi_B|^2=2|d^\nabla B|^2$, where $|d^\nabla B|$ is the norm of $d^\nabla B$ as a $T^*\Sigma$-valued two-form.
\end{proof}

We can now state the identity on the level sets. On each level set $\Sigma_r$, $r>0$, with $\mathring{\Ric}_T$ as in Lemma \ref{lem_ric_sigma_five}, let $\Phi_\Sigma=\Phi_{\mathring{\Ric}_T}$ be its Codazzi defect, defined as in Lemma \ref{lem_codazzi_int} with respect to the connection of $\Sigma_r$, and set
\begin{align}\label{eq_E_def}
	Q=\Ric_T^2-\lambda\Ric_T,\qquad \widehat W=2\lambda W_\Sigma-\frac{R_\Sigma}{2\lambda}Z_{\mathring{\Ric}_T},
\end{align}
where $Z_{\mathring{\Ric}_T}$ is defined before \eqref{eq_ZB}. The tensors $\Ric_T'$, $\Phi_\Sigma$, $Q$ and $\widehat W$ vanish on the level sets of the rigid soliton $\Sph^2(1/\sqrt\lambda)\times\R^3$, which are the products $\Sph^2(1/\sqrt\lambda)\times\Sph^2(r)$. The coefficient of $Z_{\mathring{\Ric}_T}$ in $\widehat W$ is determined as follows. In the proof of Theorem \ref{thm_identity}, the Chern--Gauss--Bonnet formula supplies the integral of $\|W_\Sigma\|^2$, while the Bochner and Codazzi identities supply the integral of $W_\Sigma(B,B)$, where $B=\mathring\Ric_T$, and $Z_B$ represents this pairing by \eqref{eq_ZB_norms}. After these substitutions, the coefficient of $|\Sigma_r|$ in the integral of $\|2\lambda W_\Sigma-\eta Z_B\|^2$ is $\frac{\lambda^2}{3}(R_\Sigma-2\lambda\eta)^2$, and it vanishes exactly when $\eta=R_\Sigma/(2\lambda)$.

\begin{theorem}\label{thm_identity}
Let $(M^5,g,f,\lambda)$ be a complete gradient shrinking Ricci soliton with constant scalar curvature $R=2\lambda$. Then, for every $r>0$,
\begin{align}\label{eq_level_intro}
	8\pi\lambda^2\int_N|\II|^2\dmu_N=\int_{\Sigma_r}\Big((R_\Sigma+2\lambda)|\Ric_T'|^2+R_\Sigma|\Phi_\Sigma|^2+\frac{R_\Sigma^2}{\lambda^2}|Q|^2+\|\widehat W\|^2\Big)\dmu_r,
\end{align}
where $R_\Sigma=2\lambda+2r^{-2}$, and $\|\widehat W\|$ is the norm of $\widehat W$ as a self-adjoint endomorphism of $\Lambda^2T\Sigma_r$, as in \eqref{eq_cgb}.
\end{theorem}

\begin{proof}
Fix $r>0$, and write
\begin{gather*}
	s=R_\Sigma=2\lambda+\frac{2}{r^2},\qquad I_0=\int_N|\II|^2\dmu_N,\qquad J_0=\int_{\Sigma_r}|\Ric_T'|^2\dmu_r.
\end{gather*}
By \eqref{eq_ricT_five}, $Q=\mathring{\Ric}_T^2-\frac{\lambda^2}{4}\Id$. All integrals below are over $\Sigma_r$, and we write $B=\mathring{\Ric}_T$ inside this proof.

By Proposition \ref{prop_weyl_norm},
\begin{align}\label{eq_W2}
	\int\|W_\Sigma\|^2=\frac{s^2|\Sigma_r|}{12}+2\pi I_0+\frac{J_0}{2\lambda^2r^2}.
\end{align}
By \eqref{eq_four_decomp} with $A=B$, $\trg B=0$ and $|B|^2=\lambda^2$, and by \eqref{eq_traceless_sigma},
\[
	\Rm_\Sigma(B,B)=W_\Sigma(B,B)-\trg(\mathring{\Ric}_\Sigma B^2)-\frac{s\lambda^2}{12},\qquad \trg(\Ric_\Sigma B^2)=\trg(\mathring{\Ric}_\Sigma B^2)+\frac{s\lambda^2}{4}.
\]
Writing $\Ric_T=B+\frac\lambda2\Id$ in an eigenbasis of $\Ric_T$, we obtain $\Rm_\Sigma(\Ric_T,\Ric_T)=\Rm_\Sigma(B,B)+\lambda\inner{\Ric_\Sigma}{B}+\frac{\lambda^2}{4}s$, and, by \eqref{eq_traceless_sigma}, $\inner{\Ric_\Sigma}{B}=\inner{\mathring{\Ric}_\Sigma}{B}=\beta_r\lambda^2=\lambda^2-\lambda r^{-2}$. Since $\nabla^\Sigma\Ric_T=\nabla^\Sigma B$ and $2\lambda^3+2\lambda^2r^{-2}=\lambda^2s$, \eqref{eq_bochner_level} with $k=2$ gives
\begin{align}\label{eq_bochner_five}
	|\nabla^\Sigma B|^2+|\Ric_T'|^2=2\Rm_\Sigma(B,B)-\frac{\lambda^2}{2}s.
\end{align}
Hence $2\int|\nabla^\Sigma B|^2=4\int\Rm_\Sigma(B,B)-2J_0-s\lambda^2|\Sigma_r|$, and Lemma \ref{lem_codazzi_int}, applied to the divergence free tensor $B$ on the closed manifold $\Sigma_r$, gives
\[
	\int|\Phi_\Sigma|^2=2\int|\nabla^\Sigma B|^2+2\int\trg(\mathring{\Ric}_\Sigma B^2)+\frac{s\lambda^2}{2}|\Sigma_r|-2\int\Rm_\Sigma(B,B).
\]
Substituting the first identity into the second,
\begin{align*}
	\int|\Phi_\Sigma|^2&=2\int\Rm_\Sigma(B,B)+2\int\trg(\mathring{\Ric}_\Sigma B^2)-2J_0-\frac{s\lambda^2}{2}|\Sigma_r|\\
	&=2\int W_\Sigma(B,B)-2J_0-\frac23\lambda^2s|\Sigma_r|,
\end{align*}
that is,
\begin{align}\label{eq_WBB}
	2\int W_\Sigma(B,B)=\int|\Phi_\Sigma|^2+2J_0+\frac23\lambda^2s|\Sigma_r|.
\end{align}
We now expand the square of $\widehat W$ using \eqref{eq_ZB_norms}, \eqref{eq_W2} and \eqref{eq_WBB}:
\begin{align*}
	\int\|\widehat W\|^2&=4\lambda^2\int\|W_\Sigma\|^2+\frac{s^2}{4\lambda^2}\int\|Z_B\|^2-2s\int W_\Sigma(B,B)\\
	&=\frac{\lambda^2s^2|\Sigma_r|}{3}+8\pi\lambda^2I_0+\frac{2J_0}{r^2}+\frac{\lambda^2s^2|\Sigma_r|}{3}-\frac{s^2}{\lambda^2}\int|Q|^2\\
	&\quad-\frac23\lambda^2s^2|\Sigma_r|-2sJ_0-s\int|\Phi_\Sigma|^2\\
	&=8\pi\lambda^2I_0-(s+2\lambda)J_0-s\int|\Phi_\Sigma|^2-\frac{s^2}{\lambda^2}\int|Q|^2,
\end{align*}
since $2s-2r^{-2}=s+2\lambda$. This is \eqref{eq_level_intro}.
\end{proof}

\begin{theorem}\label{thm_level}
Let $(M^5,g,f,\lambda)$ be a complete gradient shrinking Ricci soliton with constant scalar curvature $R=2\lambda$, and suppose that the minimum set of $f$ is totally geodesic. Then $M$ is isometric to a quotient of $\Sph^2(1/\sqrt{\lambda})\times\R^3$ by a group of order at most two.
\end{theorem}

\begin{proof}
Since $N$ is totally geodesic, the left-hand side of \eqref{eq_level_intro} vanishes for every $r>0$. All the coefficients on its right-hand side are positive, so $Q=0$ and $\Ric_T'=0$ on every level set $\Sigma_r$. Let $U=rS$. By \eqref{eq_ricT_five}, $U=\Id-\lambda^{-1}\Ric_T$, so that $Q=\lambda^2(U^2-U)$ and $U'=-\lambda^{-1}\Ric_T'$. Hence $U^2=U$ and $U'=0$, and
\[
	S'+S^2=\frac{U'}{r}+\frac{U^2-U}{r^2}=0.
\]
By the Riccati equation, $\Rm(E_1,Y,E_1,T)=0$ for $Y,T$ tangent to the level sets, so $M$ is radially flat on $M\setminus N$, and also on $N$, where $\nabla f=0$. Hence $M$ is rigid by \cite[Theorem 1.2]{pewy}. Thus the universal cover of $M$ is $\widetilde N\times\R^3$, where $\widetilde N$ is a complete simply connected surface with $\Ric_{\widetilde N}=\lambda g_{\widetilde N}$, that is, the sphere $\Sph^2(1/\sqrt\lambda)$, and $M$ is the quotient of this product by its fundamental group. Since $M$ is homotopy equivalent to $N$, which is a sphere or a real projective plane, its fundamental group has order at most two.
\end{proof}

\section{The case \texorpdfstring{$R=3\lambda$}{R = 3 lambda} in dimension five}\label{sec_three}

In this section, $(M^5,g,f,\lambda)$ is a complete gradient shrinking Ricci soliton with constant scalar curvature $R=3\lambda$ and nonconstant potential function, so that $k=3$ and $q=2$. The statement of the following theorem is that of \cite[Theorem 1.1]{loqw1}, where this case was considered by Li, Ou, Qu and Wu. Our argument is based on a cubic defect of the rescaled shape operator of the level sets and a monotonicity argument, and is independent of Theorem \ref{thm_focal}. The monotonicity argument follows the one in the proof of \cite[Corollary 1.6]{horlevel}, where a defect of a rescaled principal curvature, weighted by a power of $r$ and a Gaussian factor, is shown to vanish, which also leads to $S'+S^2=0$. The identity \eqref{eq_weyl_three}, which follows from the Chern--Gauss--Bonnet formula, also appears in \cite[Lemma 5.2 and (5.22)]{loqw1} in our notation.

\begin{theorem}\label{thm_three_five}
Let $(M^5,g,f,\lambda)$ be a complete connected gradient shrinking Ricci soliton with constant scalar curvature $R=3\lambda$. Then $M$ is a finite quotient of $\Sph^3(\sqrt{2/\lambda})\times\R^2$.
\end{theorem}

\begin{proof}
We use the notation of Section \ref{sec_level}, and we set $U=rS$, the rescaled shape operator of the level sets. By \eqref{eq_ricT} and \eqref{eq_rH},
\begin{align}\label{eq_U}
	\Ric_T=\lambda(\Id-U),\qquad \trg U=|U|^2=1,\qquad H=\frac1r.
\end{align}
Since $S=U/r$, \eqref{eq_ric_sigma} gives
\[
	\Ric_\Sigma=\Ric_T+HS+S'=\lambda(\Id-U)+\frac1rU'.
\]
Since $\trg U'=0$ and $\trg(UU')=\frac12(\trg U^2)'=0$, it follows that $R_\Sigma=3\lambda$ and
\begin{align}\label{eq_traceless_three}
	|\mathring{\Ric}_\Sigma|^2=\frac34\lambda^2+\frac{|U'|^2}{r^2}.
\end{align}
The level sets are circle bundles over the closed three-manifold $N$, so their Euler characteristic vanishes, and the Chern--Gauss--Bonnet formula \eqref{eq_cgb}, with the endomorphism norm of $W_\Sigma$ as in Section \ref{sec_level}, gives
\begin{align}\label{eq_weyl_three}
	\int_{\Sigma_r}\|W_\Sigma\|^2\dmu_r=\int_{\Sigma_r}\Big(\frac12|\mathring{\Ric}_\Sigma|^2-\frac{R_\Sigma^2}{24}\Big)\dmu_r=\frac{1}{2r^2}\int_{\Sigma_r}|U'|^2\dmu_r.
\end{align}
By Proposition \ref{prop_structure}, the volume measure of $\Sigma_r$ is $r\,\mathrm{d}\sigma$, as in Section \ref{sec_level}, so that $|\Sigma_r|=2\pi r|N|$, and the normalized measure $|\Sigma_r|^{-1}\mathrm{d}\mu_r$ is the fixed probability measure $(2\pi|N|)^{-1}\mathrm{d}\sigma$ on the unit normal bundle of $N$. We write
\[
	\langle u\rangle_r=\frac{1}{|\Sigma_r|}\int_{\Sigma_r}u\dmu_r,
\]
so that $\frac{\mathrm{d}}{\mathrm{d}r}\langle u\rangle_r=\langle u'\rangle_r$. Then \eqref{eq_weyl_three} reads $\langle|U'|^2\rangle_r=2r^2\langle\|W_\Sigma\|^2\rangle_r$.

By Lemma \ref{lem_bochner_level} with $k=3$ and \eqref{eq_U}, $\lambda^2\big(|\nabla^\Sigma U|^2+|U'|^2\big)=2\Rm_\Sigma(\Ric_T,\Ric_T)-6\lambda^3$, and
\[
	\Rm_\Sigma(\Ric_T,\Ric_T)=\lambda^2\big(R_\Sigma-2\inner{\Ric_\Sigma}{U}+\Rm_\Sigma(U,U)\big)=\lambda^2\big(3\lambda+\Rm_\Sigma(U,U)\big),
\]
since $\inner{\Ric_\Sigma}{U}=\lambda(\trg U-\trg U^2)+r^{-1}\trg(UU')=0$. Hence
\begin{align}\label{eq_U_bochner}
	|U'|^2+|\nabla^\Sigma U|^2=2\Rm_\Sigma(U,U).
\end{align}
By \eqref{eq_four_decomp} and \eqref{eq_U},
\[
	\Rm_\Sigma(U,U)=W_\Sigma(U,U)+\trg\big(\mathring{\Ric}_\Sigma U(\Id-U)\big)=W_\Sigma(U,U)-\inner{\Ric_\Sigma}{U^2}.
\]
Define the cubic defect
\[
	v=1-\trg U^3,\qquad \theta(r)=\langle v\rangle_r.
\]
Since $\Ric_T=\lambda(\Id-U)$ and $\Ric(E_1,\cdot)=0$, we have $\lambda^3v=\trg(\Ric^3)-\lambda|\Ric|^2$, and, with $\lambda=\frac12$, this scalar also appears in \cite[Section 6]{loqw1}. Every eigenvalue $u_i$ of $U$ belongs to $[-1,1]$, since $|U|^2=1$, so that
\begin{align}\label{eq_v_bounds}
	0\leq v=\sum_iu_i^2(1-u_i)\leq2,\qquad 0\leq\theta\leq2.
\end{align}
Moreover,
\[
	\inner{\Ric_\Sigma}{U^2}=\lambda(\trg U^2-\trg U^3)+\frac1r\trg(U^2U')=\lambda v-\frac{v'}{3r},
\]
and averaging \eqref{eq_U_bochner} gives
\begin{align}\label{eq_V_ode}
	\theta'-3\lambda r\theta=\frac{3r}{2}\big\langle|\nabla^\Sigma U|^2+|U'|^2-2W_\Sigma(U,U)\big\rangle_r.
\end{align}

Set $B=U-\frac14\Id$, which is trace free. Since $\trg U=\trg U^2=1$, we have
\[
	|B|^2=\frac34,\qquad \trg B^4=\trg U^4-\trg U^3+\frac{21}{64},
\]
and \eqref{eq_ZB_norms} gives
\begin{align}\label{eq_v_identity}
	v=1-\trg U^3=|U^2-U|^2+\frac14\|Z_B\|^2.
\end{align}
In particular, $v\geq0$, and $v=0$ if and only if $U^2=U$, since $U^2=U$ gives $\trg U^3=\trg U=1$. Since $W_\Sigma$ is trace free, $W_\Sigma(U,U)=W_\Sigma(B,B)=\inner{W_\Sigma}{Z_B}_{\mathrm{HS}}$ by \eqref{eq_ZB_norms}. We replace $\langle|U'|^2\rangle_r$ by $2r^2\langle\|W_\Sigma\|^2\rangle_r$ in \eqref{eq_V_ode}, complete the square, and use \eqref{eq_v_identity} in the form $\frac{3}{4r}\|Z_B\|^2=\frac3r\big(v-|U^2-U|^2\big)$. This gives
\begin{align}\label{eq_W_average}
	\theta'-\Big(3\lambda r-\frac3r\Big)\theta=\frac{3r}{2}\Big\langle|\nabla^\Sigma U|^2+2r^2\Big\|W_\Sigma-\frac{Z_B}{2r^2}\Big\|^2\Big\rangle_r+\frac3r\big\langle|U^2-U|^2\big\rangle_r.
\end{align}
In particular, $\theta'\geq(3\lambda r-\frac3r)\theta$, so that $z(r)=r^3e^{-3\lambda r^2/2}\theta(r)$ is nonnegative and nondecreasing on $(0,\infty)$. By \eqref{eq_v_bounds}, $z(r)\to0$ as $r\to\infty$. Hence $z=0$, and $\theta=0$ identically.

Since $v\geq0$ is continuous, it vanishes pointwise, so that $U^2=U$ by \eqref{eq_v_identity}, and $U$ is a rank one orthogonal projection, since $\trg U=1$. Hence $\Rm_\Sigma(U,U)=0$, and \eqref{eq_U_bochner} gives $U'=0$ and $\nabla^\Sigma U=0$. Since $S=U/r$,
\[
	S'+S^2=\frac{U'}{r}+\frac{U^2-U}{r^2}=0,
\]
and the Riccati equation of Section \ref{sec_level} gives $\Rm(E_1,Y,E_1,T)=0$ for $Y,T$ tangent to the level sets. Hence $M$ is radially flat on $M\setminus N$, and also on $N$, where $\nabla f=0$, and $M$ is rigid by \cite[Theorem 1.2]{pewy}. Its universal cover is $\widetilde N\times\R^2$, where $\widetilde N$ is a complete simply connected three-manifold with $\Ric_{\widetilde N}=\lambda g_{\widetilde N}$, that is, the round sphere of sectional curvature $\lambda/2$, and the quotient is finite, since the deck group acts freely and properly discontinuously on the compact set $\widetilde N\times\{0\}$, which it preserves.
\end{proof}

\begin{proof}[Proof of Theorem \ref{thm_five}]
By \cite[Theorems 1 and 10]{fr16}, $R\in\{0,2\lambda,3\lambda,4\lambda,5\lambda\}$. If $R=0$ or $R=5\lambda$, $M$ is the Gaussian soliton or a compact Einstein manifold \cite[Propositions 3.1 and 3.3]{pewy}. If $R=4\lambda$, $M$ is rigid by \cite[Theorem 10]{fr16}, and its universal cover is $\widetilde N^4\times\R$ with $\widetilde N^4$ a compact simply connected Einstein manifold with $\Ric_{\widetilde N}=\lambda g_{\widetilde N}$. If $R=3\lambda$, $M$ is a finite quotient of $\Sph^3(\sqrt{2/\lambda})\times\R^2$ by Theorem \ref{thm_three_five}. Suppose now that $R=2\lambda$. Then $f$ is not constant, and, by Theorem \ref{thm_focal}, the minimum set $N$ of $f$ is totally geodesic. By Theorem \ref{thm_level}, $M$ is a quotient of $\Sph^2(1/\sqrt\lambda)\times\R^3$ by a group of order at most two.
\end{proof}

\section{Four-dimensional focal sets}\label{sec_four}

In this section, $k=4$, so that $N$ is a closed minimal four-manifold of codimension $q=n-4$, and we prove Theorems \ref{thm_four} and \ref{thm_harmonic_weyl}. We first explain why the method of Section \ref{sec_lowdim} requires an additional hypothesis.

\begin{remark}\label{rem_four}
For $k\in\{2,3\}$, the intrinsic curvature decomposition of $N$ contains no Weyl term. For $k\geq4$, the identity of Theorem \ref{thm_general} contains the pairing $W^N(\II,\II)$, which has no prescribed sign. Even among rigid solitons, the focal set can be any closed Einstein manifold $N^k$ with $\Ric_N=\lambda g_N$, since $N\times\R^{n-k}$ is such a soliton, and constant scalar curvature does not force $W^N=0$, as the normalized products with focal set $\Sph^2\times\Sph^2$ or $\mathbb{CP}^2$ show. Hence a Weyl restriction which fails on a rigid model cannot follow from constant scalar curvature alone. This does not exclude control of the pairing $W^N(\II,\II)$, which vanishes on every rigid model. Two natural intrinsic hypotheses on the metric induced on $N$ are $\delta_NW^N=0$ and $B_N=0$, where $B_N$ is the Bach tensor of $(N,g_N)$, and both hold when $g_N$ is Einstein. Further rigidity questions under harmonic Weyl curvature or Bach-flatness of the metric induced on $N$ are being investigated by the author in work in preparation.
\end{remark}

We first prove Theorem \ref{thm_four}.

\begin{proof}[Proof of Theorem \ref{thm_four}]
Here $k=4$ and $q=n-4$. Write $w=W^N(\II,\II)$. With $k=4$, \eqref{eq_simons_weyl} reads $\trg(\Ric_NT)+\omega_{\II}=-w+2\inner{\mathring{\Ric}_N}{\mathring{T}}+\frac43\lambda t+\frac16t^2$. We subtract three times the vanishing integral $\int_NI\dmu_N$ of Lemma \ref{lem_simons} from the integral of \eqref{eq_general_weyl} with $k=4$, which also vanishes by Theorem \ref{thm_general}. Using $C_D=|D\II|^2-\frac12|A|^2$, $O_{\II}=\frac12\inner{K}{V}+C$, $|\mathring{T}|^2=|T|^2-\frac14t^2$ and
\begin{gather*}
	2\Big|V-\frac34K\Big|^2+\frac32\inner{K}{V}=2\Big|V-\frac38K\Big|^2+\frac{27}{16}C,\\
	6|\mathring{\Ric}_N|^2+2\inner{\mathring{\Ric}_N}{\mathring{T}}=6\Big|\mathring{\Ric}_N+\frac16\mathring{T}\Big|^2-\frac16|\mathring{T}|^2,
\end{gather*}
we obtain
\begin{align}\label{eq_four_identity}
\begin{split}
	0=\int_N\Big(&|D\II|^2+\frac92|A|^2-3|d|^2+4\Big|L^s+\frac12P\Big|^2+2\Big|V-\frac38K\Big|^2+2|Y|^2+\frac{1}{\lambda^2}|F|^2\\
	&+6\Big|\mathring{\Ric}_N+\frac16\mathring{T}\Big|^2+\frac23\lambda t-11w+\frac{7}{16}C+\frac72c-\frac14t^2-\frac16|T|^2\Big)\dmu_N.
\end{split}
\end{align}
By \eqref{eq_d_bound} with $k=4$, $\frac92|A|^2-3|d|^2\geq0$. The shape operators lie in the nine-dimensional space of trace free symmetric endomorphisms of $T_pN$, so $G$ is a positive semidefinite Gram matrix of rank at most $\min\{q,9\}$ and trace $t$, whence $c\geq\frac19t^2$, and $c\geq\frac12t^2$ if $q=2$. If $a_1,\ldots,a_4$ are the eigenvalues of a trace free symmetric endomorphism $S$, then $a_i^2=(\sum_{j\neq i}a_j)^2\leq3\sum_{j\neq i}a_j^2$, so that $S^2\leq\frac34|S|^2\Id$. Hence $0\leq T\leq\frac34t\Id$ and $|T|^2\leq\frac34t\trg T=\frac34t^2$, and
\[
	\frac{7}{16}C+\frac72c-\frac14t^2-\frac16|T|^2\geq\Big(\frac{7}{18}-\frac14-\frac18\Big)t^2=\frac{1}{72}t^2,
\]
with $\frac{7}{18}$ replaced by $\frac74$, and $\frac{1}{72}$ by $\frac{11}{8}$, if $q=2$. Discarding the nonnegative terms in \eqref{eq_four_identity} gives \eqref{eq_four_estimate}. If \eqref{eq_four_pinching} holds, then \eqref{eq_four_estimate} gives $\int_N(|D\II|^2+\frac{1}{72}|\II|^4)\dmu_N\leq0$, so $\II=0$, and $M$ is rigid by \cite[Corollary 1.1]{lloqwz}.
\end{proof}

\begin{remark}\label{rem_multiplier}
For $k\geq3$, the Simons-type identity of Lemma \ref{lem_simons} supplies a second vanishing integral, with integrand $I$. Thus we may replace $J$ by $J-sI$ for a constant $s$, and this explains the multiples of $I$ in the proofs of Theorems \ref{thm_three} and \ref{thm_four}. In the notation of Section \ref{sec_bochner}, the derivative terms of $J-sI$ and its coefficient of $\lambda t$ are, respectively,
\[
 (4-s)|D\II|^2+\Big(3+\frac{s}{2}\Big)|A|^2+(s-6)|d|^2,
 \qquad \frac{14-sk}{k-1}.
\]
When $k=3$, taking $s=4$ eliminates $|D\II|^2$ and leaves $\lambda t$ with coefficient one, and completing the two remaining quadratic expressions gives the squares with shifts $\frac14K$ and $\frac56\mathring T$ in Theorem \ref{thm_three}. When $k=4$, the same choice would make the coefficient of $\lambda t$ negative. For $3\leq s\leq\frac72$, the bound \eqref{eq_d_bound} makes the derivative terms nonnegative, and the estimates on $c$ and $T$ used in the proof of Theorem \ref{thm_four} bound the remaining quartic terms below by
\[
 \left(\frac1{72}-\frac{(s-3)^2}{12}\right)t^2.
\]
This bound is largest at $s=3$, which gives \eqref{eq_four_identity}. The contraction $W^N(\II,\II)$ remains, with coefficient $s-14$, and accounts for the hypothesis of Theorem \ref{thm_four}; Theorem \ref{thm_harmonic_weyl} follows by applying the Bochner formula for the intrinsic Weyl tensor.
\end{remark}

With the conventions of Theorem \ref{thm_harmonic_weyl}, if $\delta_NW^N=0$, then
\begin{align}\label{eq_weyl_bochner}
	\frac12\Delta_N\|W^N\|^2=\|\nabla^NW^N\|^2+\frac{R_N}{2}\|W^N\|^2-6\trg_{\Lambda^2}\big((W^N)^3\big),
\end{align}
see \cite[Corollary 4.2]{catmas}.

\begin{proof}[Proof of Theorem \ref{thm_harmonic_weyl}]
Here $k=4$. By Lemma \ref{lem_hess_two}, $R_N=4\lambda+\frac12t$. Since \eqref{eq_weyl_bochner} is a pointwise identity, integrating it over $N$ gives
\[
	6\int_N\trg_{\Lambda^2}\big((W^N)^3\big)\dmu_N-\int_N\|\nabla^NW^N\|^2\dmu_N-2\lambda\int_N\|W^N\|^2\dmu_N=\frac14\int_Nt\|W^N\|^2\dmu_N.
\]
By \eqref{eq_weyl_condition}, $t\|W^N\|^2=0$ at every point of $N$. Hence, at every point, $\II=0$ or $W^N=0$, so that $W^N(\II,\II)=0$, and Theorem \ref{thm_four} applies.
\end{proof}

Every rigid soliton with $R=4\lambda$ satisfies the additional curvature conditions of Theorems \ref{thm_four} and \ref{thm_harmonic_weyl}: its focal set is totally geodesic and Einstein with $\Ric_N=\lambda g_N$, so that $\delta_NW^N=0$, and \eqref{eq_weyl_bochner} with $R_N=4\lambda$ shows that \eqref{eq_weyl_condition} holds with equality.

\appendix

\section{The connection terms}\label{app_conn}

We prove the formulas \eqref{eq_Omega} for $\Omega_3$ and $\Omega_1$, which are used in the proof of Lemma \ref{lem_X}. Indices are as in Section \ref{sec_prelim}. We begin with $\Omega_3$. In the frames of the proof of Lemma \ref{lem_hessian}, we have $\nabla_{e_j}e_i=h^\beta_{ji}e_\beta$ and $\nabla_{e_j}e_\alpha=-h^\alpha_{jm}e_m$ at $p$, so that, at $p$,
\begin{align*}
	(\nabla_{e_j}\Rm)(e_i,e_\alpha,e_i,e_\ell)={}&e_j\big(R_{i\alpha i\ell}\big)-\Rm(\nabla_{e_j}e_i,e_\alpha,e_i,e_\ell)-\Rm(e_i,\nabla_{e_j}e_\alpha,e_i,e_\ell)\\
	&-\Rm(e_i,e_\alpha,\nabla_{e_j}e_i,e_\ell)-\Rm(e_i,e_\alpha,e_i,\nabla_{e_j}e_\ell)\\
	={}&e_j\big(R_{i\alpha i\ell}\big)-h^\beta_{ji}R_{\beta\alpha i\ell}+h^\alpha_{jm}R_{imi\ell}-h^\beta_{ji}R_{i\alpha\beta\ell}-h^\beta_{j\ell}R_{i\alpha i\beta}.
\end{align*}
By Lemma \ref{lem_first_order}(ii), $\sum_iR_{i\alpha i\ell}=-d^\alpha_\ell$ at every point of $N$, so $\sum_ie_j(R_{i\alpha i\ell})=-(D_{e_j}d^\alpha)_\ell$ at $p$. Multiplying by $h^\alpha_{j\ell}$ and summing over $i,j,\ell,\alpha$, we obtain
\[
	\Omega_3=\tau_3-\sum h^\beta_{ji}h^\alpha_{j\ell}R_{\beta\alpha i\ell}+\sum h^\alpha_{jm}h^\alpha_{j\ell}R_{imi\ell}-\sum h^\beta_{ji}h^\alpha_{j\ell}R_{i\alpha\beta\ell}-\sum h^\beta_{j\ell}h^\alpha_{j\ell}R_{i\alpha i\beta}.
\]
We use
\[
	\sum_jh^\beta_{ij}h^\alpha_{j\ell}=(S_\beta S_\alpha)_{i\ell}=\frac12(P_{\beta\alpha}+K_{\beta\alpha})_{i\ell},
\]
where $P_{\beta\alpha}$ is symmetric and $V_{\beta\alpha i\ell}$ and $K_{\beta\alpha}$ are skew in $(i,\ell)$. The four sums are
\begin{align*}
	\sum h^\beta_{ji}h^\alpha_{j\ell}R_{\beta\alpha i\ell}&=\sum(S_\beta S_\alpha)_{i\ell}V_{\beta\alpha i\ell}=\frac12\inner{K}{V},\\
	\sum h^\alpha_{jm}h^\alpha_{j\ell}R_{imi\ell}&=\sum(S_\alpha^2)_{m\ell}\rho_{m\ell}=\inner{T}{\rho},\\
	\sum h^\beta_{ji}h^\alpha_{j\ell}R_{i\alpha\beta\ell}&=-\sum(S_\beta S_\alpha)_{i\ell}L_{\beta\alpha\ell i}=-\frac12\inner{P}{L^t}-\frac12\inner{K}{L^t}=-\frac12\inner{P}{L^s}+\frac14\inner{K}{V},\\
	\sum h^\beta_{j\ell}h^\alpha_{j\ell}R_{i\alpha i\beta}&=-\frac32\sum G_{\alpha\beta}G_{\alpha\beta}=-\frac32c,
\end{align*}
where we used Lemma \ref{lem_hess_two}(i) and
\begin{gather*}
	R_{i\alpha\beta\ell}=R_{\beta\ell i\alpha}=-R_{\beta\ell\alpha i}=-L_{\beta\alpha\ell i},\\
	\inner{P}{L^t}=\inner{P}{L^s},\qquad \inner{K}{L^t}=\inner{K}{L^s-\tfrac12V}=-\tfrac12\inner{K}{V}.
\end{gather*}
Therefore
\[
	\Omega_3=\tau_3-\frac12\inner{K}{V}+\inner{T}{\rho}+\frac12\inner{P}{L^s}-\frac14\inner{K}{V}+\frac32c=\tau_3+\inner{T}{\rho}+\frac12\inner{P}{L^s}-\frac34\inner{K}{V}+\frac32c,
\]
which is the first formula in \eqref{eq_Omega}. We now compute $\Omega_1$ in the same way. At $p$,
\[
	(\nabla_{e_i}\Rm)(e_j,e_\alpha,e_i,e_\ell)=e_i\big(R_{j\alpha i\ell}\big)-h^\beta_{ij}R_{\beta\alpha i\ell}+h^\alpha_{im}R_{jmi\ell}-h^\beta_{ii}R_{j\alpha\beta\ell}-h^\beta_{i\ell}R_{j\alpha i\beta}.
\]
By Lemma \ref{lem_first_order}(ii), $R_{j\alpha i\ell}=A^\alpha_{\ell ij}$ at every point of $N$, so $e_i(R_{j\alpha i\ell})=(D_{e_i}A^\alpha)_{\ell ij}$ at $p$. Multiplying by $h^\alpha_{j\ell}$ and summing, we obtain
\[
	\Omega_1=\tau_1-\sum h^\beta_{ij}h^\alpha_{j\ell}R_{\beta\alpha i\ell}+\sum h^\alpha_{im}h^\alpha_{j\ell}R_{jmi\ell}-\sum_ih^\beta_{ii}\sum h^\alpha_{j\ell}R_{j\alpha\beta\ell}-\sum h^\beta_{i\ell}h^\alpha_{j\ell}R_{j\alpha i\beta},
\]
and the four sums are
\begin{align*}
	\sum h^\beta_{ij}h^\alpha_{j\ell}R_{\beta\alpha i\ell}&=\frac12\inner{K}{V},\\
	\sum h^\alpha_{im}h^\alpha_{j\ell}R_{jmi\ell}&=-\sum h^\alpha_{im}h^\alpha_{j\ell}R_{jm\ell i}=-\omega,\\
	\sum_ih^\beta_{ii}&=0,\\
	\sum h^\beta_{i\ell}h^\alpha_{j\ell}R_{j\alpha i\beta}&=\sum(S_\alpha S_\beta)_{ji}L_{\alpha\beta ji}=\frac12\inner{P}{L}+\frac12\inner{K}{L}=\frac12\inner{P}{L^s}+\frac14\inner{K}{V},
\end{align*}
where the first sum is the one computed above, the second uses the skew symmetry of $R$ in its last two entries and the definition $\omega=\sum R_{ij\ell m}h^\alpha_{i\ell}h^\alpha_{jm}$ with $(i,j,\ell,m)$ renamed $(j,m,\ell,i)$, the third is the minimality of $N$ (Proposition \ref{prop_structure}), and the fourth uses $R_{j\alpha i\beta}=L_{\alpha\beta ji}$ and $\inner{K}{L}=\frac12\inner{K}{L-L^t}=\frac12\inner{K}{V}$. Therefore
\[
	\Omega_1=\tau_1-\frac12\inner{K}{V}-\omega-\frac12\inner{P}{L^s}-\frac14\inner{K}{V}=\tau_1-\omega-\frac12\inner{P}{L^s}-\frac34\inner{K}{V}.
\]

\section{Notation}\label{app_notation}

The following table lists the quantities attached to the minimum set $N$ in Sections \ref{sec_prelim} to \ref{sec_lowdim} and \ref{sec_four}, with the place of their definition.

\begin{center}
\small
\renewcommand{\arraystretch}{1.3}
\begin{tabular}{@{}l@{\quad}l@{\qquad}l@{\quad}l@{}}
	$E$, $P_E$, $P_{TN}$ & Proposition \ref{prop_structure} & $\Theta$, $\Omega$ & \eqref{eq_bochner_laplacian}\\
	$q$ & \eqref{eq_laplacian_f} & $\omega$, $A$, $\tau_1$, $\tau_3$, $C_D$ & before Lemma \ref{lem_M}\\
	$\II$, $h^\alpha_{ij}$, $S_\alpha$, $D$, $R^\perp$ & Section \ref{sec_prelim} & $\Omega_1$, $\Omega_3$ & proof of Lemma \ref{lem_X}\\
	$t$, $G$, $c$, $T$, $P$, $K$, $\sigma$, $d$ & Section \ref{sec_prelim} & $a$, $C$ & \eqref{eq_algebraic_norms}\\
	$L$, $V$, $Y$, $L^s$, $L^t$, $\rho$, $\tilde\rho$ & Section \ref{sec_prelim} & $J$ & Theorem \ref{thm_general}\\
	$\Lambda$, $h^\alpha_{ij,\ell}$ & Section \ref{sec_hessian} & $\mathring{T}$, $\mathring{\Ric}_N$, $\mathrm{Sch}^N$, $W^N(\II,\II)$ & after Theorem \ref{thm_general}\\
	$F$ & Lemma \ref{lem_hessian} & $O_{\II}$, $\omega_{\II}$ & before Lemma \ref{lem_simons}\\
	& & $I$ & Lemma \ref{lem_simons}
\end{tabular}
\end{center}

\end{document}